\documentclass[a4paper,11pt]{article}

\usepackage{microtype}
\usepackage[]{algorithm2e}
\usepackage{chngpage}

\usepackage[T1]{fontenc}
\usepackage[utf8]{inputenc}
\usepackage[doublespacing]{setspace}

\usepackage{lineno}

\usepackage{amsfonts}
\usepackage{amsmath}
\numberwithin{equation}{section}

\usepackage{amssymb}
\usepackage{mathtools}
\usepackage[binary-units=true]{siunitx}

\usepackage{bm}
\usepackage{breqn}

\usepackage{amsthm}
\newtheorem{theorem}{Theorem} 
\newtheorem{lemma}[theorem]{Lemma}
\newtheorem{assumption}{Assumption}

\usepackage{caption}
\usepackage{subcaption}
\usepackage{diagbox}

\usepackage{txfonts}
\usepackage{sansmath}

\usepackage{geometry}
\newgeometry{vmargin={25mm}, hmargin={25mm,25mm}}

\usepackage{multicol}
\usepackage{booktabs}
\usepackage{graphicx}
\usepackage{float}

\makeatletter
\newcommand\tabfill[1]{%
\dimen@\linewidth%
\advance\dimen@\@totalleftmargin%
\advance\dimen@-\dimen\@curtab%
\parbox[t]\dimen@{#1\ifhmode\strut\fi}%
}

\usepackage{dcolumn}
\newcolumntype{d}[1]{D{.}{.}{#1}}

\usepackage{tikz}
\usetikzlibrary{shapes,arrows}
\usetikzlibrary{matrix}
\usepackage{pgfplots}
\usepackage{pgfplotstable}
								
\usepackage{color}
\definecolor{thered}{rgb}{0.65,0.04,0.07}
\definecolor{thegreen}{rgb}{0.06,0.44,0.08}
\definecolor{theblue}{rgb}{0.02,0.2,0.68}
\definecolor{clr1}{HTML}{bdd7e7}
\definecolor{clr2}{HTML}{6baed6}
\definecolor{clr3}{HTML}{3182bd}
\definecolor{clr4}{HTML}{08519c}

\newtheorem{remark}{Remark}

\renewcommand{\@algocf@capt@plain}{above}% formerly {bottom}

\usepackage[unicode=true,
            bookmarks=true,
            bookmarksnumbered=false,
            bookmarksopen=false,
            breaklinks=true,
            pdfborder={0 0 1},
            backref=false,
            colorlinks=true,
            linkcolor=black,
            urlcolor=theblue,
            citecolor=theblue,
            hyperfootnotes=false]
           {hyperref}
\hypersetup{pdftitle={},
            pdfauthor={B. C. Vermeire}}

\begin{document}
%\linenumbers

\title{On the Strong Stability Preserving Property of Runge-Kutta Methods for Hyperbolic Problems}

\author{Mohammad R. Najafian and Brian C. Vermeire\\\\
  \textit{Department of Mechanical, Industrial, and Aerospace Engineering}\\
	\textit{Concordia University}\\
	\textit{Montreal, QC, Canada}}
\maketitle

\begin{abstract}

Strong Stability Preserving (SSP) time integration schemes maintain stability of the forward Euler method for any initial value problem. However, only a small subset of Runge-Kutta (RK) methods are SSP, and many efficient high-order time integration schemes do not formally belong to this class. In this work, we introduce a mathematical strategy to analyze the nonlinear stability of RK schemes that may not necessarily belong to the SSP class. With this approach, we mathematically demonstrate that there are time integration schemes outside the class of SSP schemes that can maintain entropy stability and positivity of density and pressure for the Lax-Friedrichs discretization, and Total Variation Diminishing stability for the first-order upwind and the second-order MUSCL schemes. As a result, for these problems, a broader range of RK methods, including the classical fourth-order, four-stage RK scheme, can be used while the numerical integration remains stable. Numerical experiments confirm these theoretical findings, and additional experiments demonstrate similar observations for a wider class of space discretizatins.

\end{abstract}

\section{Introduction}

Hyperbolic conservation laws arise in a wide range of physical applications, notably in fluid dynamics. Solutions of nonlinear hyperbolic partial differential equations (PDEs) may involve sharp gradients and discontinuities, which may cause stability issues when numerically solving these PDEs. To improve robustness, numerical methods can be designed to satisfy stability in certain convex functionals, coming from the physics of the problem. For instance, for systems of hyperbolic conservation laws, the physically meaningful solution should remain stable with respect to a nontrivial entropy function \cite{leflochHyperbolicSystemsConservation2002}. Additionally, physical quantities such as density and internal energy in compressible Euler and Navier-Stokes equations must remain positive \cite{zhangMaximumprinciplesatisfyingPositivitypreservingHighorder2011}. For scalar conservation laws, the physical solution is entropy stable, total variation diminishing, and satisfies a strict maximum principle \cite{zhangMaximumprinciplesatisfyingHighOrder2010}.

To fully discretize hyperbolic conservation laws, the method of lines is a common approach. With this method, the spatial domain is discretized first to obtain time-dependent ordinary differential equations (ODEs)
\begin{subequations}\label{IVP}
  \begin{align}
    & \frac{d}{dt}\underline{q} \left(t \right) = R \left( \underline{q} \left( t \right) \right) ,\\
    & \underline{q} \left(0 \right) = \underline{q}^0 ,
  \end{align}
\end{subequations} 
where $\underline{q} \in \mathbb{R}^N$ is the vector of unknowns, and $R\left( \underline{q}\right)$ is the right-hand side (RHS) function resulting from the spatial discretization. This ODE system is then integrated numerically in time to obtain a fully discrete solution in space and time. 

To maintain nonlinear stability properties of the continues problem, spatial discretization schemes that preserve stability when coupled with the forward Euler time integration scheme are developed \cite{rueda-ramirezSubcellLimitingStrategies2022, kuzminBoundpreservingFluxLimiting2022, linPositivityPreservingStrategy2023a}. Let $G\left( \underline{q} \right)$ represent a norm, semi-norm, or convex functional, and assume that the spatial discretization method results in a strongly stable solution in $G$ when the ODE system \ref{IVP} is integrated with the first-order forward Euler method with a time step restricted by $\Delta t_{FE}$
\begin{subequations}
	\begin{align}
	& \underline{q}^{FE}= \underline{q}^n + \Delta t R \left( \underline{q}^n \right), \\
	& G \left( \underline{q}^{FE} \right) \leq G \left( \underline{q}^n \right), \quad \mbox{for} \quad \Delta t \leq \Delta t_{FE},
	\end{align}
\end{subequations} 
where $\underline{q}^n$ is an approximation of solution at $t=t_n$, and $\Delta t_{FE}$ depends only on the RHS function.

The challenge is that even if the first-order forward Euler method provides a strongly stable solution in $G$, such stability is not in general guaranteed when higher order time integration schemes such as Runge-Kutta (RK) methods are employed \cite{gottliebStrongStabilityPreservingHighOrder2001b}. With a standard $s$-stage RK method, the next step solution, $\underline{q}^{RK}$, is approximated as  
\begin{subequations}\label{RK}
\begin{align}
& \underline{q}(t_n + \Delta t) \approx  \underline{q}^{RK} = \underline{q}^{n} + \Delta t \sum_{j=1}^s b_j \underline{R}^j, \\
&\underline{R}^i= R\left( \underline{q}^i \right), \quad \underline{q}^i = q^n + \Delta t \sum_{j=1}^{s} a_{ij} \underline{R}^j, \quad c_i= \sum_{j=1}^s a_{ij},
\end{align}
\end{subequations}
where $a_{ij}$ and $b_j$ are the coefficients of the selected RK method, and $\underline{q}^i$ are RK stage solutions. We consider in this work explicit RK schemes, where $a_{ij}=0$ for $j\geq i$. For hyperbolic PDEs, it is often desirable to maintain strong stability of the forward Euler method with intermediate stage solutions and the step solution, i.e., $\left\{G\left(\underline{q}^i \right), G\left( \underline{q}^{RK} \right)\right\} \leq G\left( \underline{q}^n \right) $. 

To study stability of RK methods when the forward Euler solution is strongly stable, it is common to write RK methods in the so-called Shu-Osher form \cite{shuEfficientImplementationEssentially1988a}
\begin{subequations}\label{Shu-Osher}
\begin{align}
\underline{q}^1 & = \underline{q}^n, \\
\underline{q}^i &= \sum_{j=1}^{i-1} \left( \alpha_{ij} \underline{q}^j + \Delta t \beta_{ij} \underline{R}^j \right), \quad i= 2, ..., s+1, \\
\underline{q}^{RK} &= \underline{q}^{s+1},
\end{align}
\end{subequations}
where consistency requires $\sum_{j=1}^{i-1}\alpha_{ij}=1$ \cite{gottliebTotalVariationDiminishing1998}. When all $\alpha_{ij}$ and $\beta_{ij}$ coefficients are non-negative, and $\alpha_{ij}$ is zero only if $\beta_{ij}$ is zero, each $\underline{q}^i$ in Eq. \ref{Shu-Osher} can be interpreted as a convex combination of forward Euler steps with a modified step size
\begin{equation}
\underline{q}^i = \sum_{j=1}^{i-1} \alpha_{ij} \left(  \underline{q}^j + \Delta t \frac{\beta_{ij}}{\alpha_{ij}}  \underline{R}^j \right).
\end{equation}
In this way, it is well-known that a sufficient condition to preserve the strong stability of forward Euler method is \cite{gottliebTotalVariationDiminishing1998}
\begin{equation}
\Delta t \leq \mathcal{C} \Delta t_{FE},  \quad \mathcal{C} = \min_{i,j} \frac{\alpha_{ij}}{\beta_{ij}},
\end{equation}
where $\alpha_{ij}/\beta_{ij}$ is interpreted as infinite when $\beta_{ij}=0$. The maximum value of $\mathcal{C}$ for a given RK method, known as the strong stability preserving (SSP) coefficient and denoted by $c^{ssp}$, can be determined by Butcher tableau coefficients of that method \cite[Theorem 3.2.]{gottliebStrongStabilityPreserving2011} \cite{ketchesonPracticalImportanceSSP2005, ferracinaStepsizeRestrictionsTotalVariationDiminishing2004}. If the SSP coefficient of an RK method is greater than zero, the scheme belongs to a family of time integration schemes known as SSP methods. An SSP method guarantees maintaining strong stability of the forward Euler method for any ODE system, initial condition, and convex stability function, provided that $\Delta t \leq c^{ssp}\Delta t_{FE}$. The necessity of this step size restriction to preserve strong stability of the forward Euler method has been shown in practice using some designed initial value problems \cite[Theorem 3.3., Example 3.7.]{gottliebStrongStabilityPreserving2011}. However, there remains a lack of connection between the theoretical proofs showing the necessity of the SSP time step size limit and practical hyperbolic problems of interest \cite{ketchesonPracticalImportanceSSP2005}.

SSP time step restriction significantly limits utilization of many RK schemes, as only a small subclass of RK methods have positive SSP coefficients. For instance, any fourth-order, four-stage RK method, including the classical RK44 scheme, has an SSP coefficient of zero and therefore does not belong to the class of SSP schemes \cite{gottliebTotalVariationDiminishing1998}. With non-SSP schemes, if no additional assumptions are made, maintaining stability of the forward Euler method is not guaranteed for any positive step size \cite[Theorem 3.3.]{gottliebStrongStabilityPreserving2011}. Having an SSP coefficient of zero has largely restricted the use of such non-SSP high-order RK methods in hyperbolic problems where maintaining stability of forward Euler method is essential.

On the other hand, there are example cases in the literature showing that, in practice, well-designed RK methods, even with an SSP coefficient of zero, maintained stability of the forward Euler method with nonzero step sizes when applied to hyperbolic conservation laws. Ketcheson and Robinson \cite{ketchesonPracticalImportanceSSP2005} showed that the midpoint method, while not being SSP, preserved Total Variation (TV) stability for Burgers' equation. The same study also reported that certain non-SSP schemes performed comparable to SSP ones when applied to Euler equations with relatively small step sizes. Good performance of some non-SSP schemes in terms of TV stability or positivity preservation for Burgers' and Euler equations are reported in \cite{gottliebHighOrderStrong2009a}. Preserving positivity with time steps larger than theoretical SSP step size limit is reported in \cite{higuerasPositivityPropertiesClassical2010, horvathPositivityStepSize2005}.

Despite the good performance of some non-SSP schemes in practical applications, it is still considered safer to avoid non-SSP schemes whenever maintaining stability of the forward Euler method is necessary \cite{ketchesonPracticalImportanceSSP2005, gottliebStrongStabilityPreservingHighOrder2001b}. The necessity of respecting the SSP time step limitation for hyperbolic problems remains an open question, as the current SSP theory does not provide sufficient information about the stability properties of non-SSP schemes for specific problems or stability functions. To study the stability of general RK methods when applied to hyperbolic problems, we introduce in the following section a novel representation of RK methods, distinct from the Shu-Osher form (Eq. \ref{Shu-Osher}). Then, with this modified representation we show that a broad class of RK schemes, while not necessarily being SSP, can maintain entropy stability of the forward Euler method with nonzero step size limits. Moreover, for certain problems, we demonstrate that these RK schemes also maintain positivity of density, positivity of internal energy, and the total variation diminishing property under nonzero time-step restrictions. Several related numerical examples are presented in Section \ref{Section_examples}.

\section{SSP property of general RK methods}

To study the SSP property of general RK methods, we introduce a modified representation of RK methods
\begin{subequations}\label{Modified_RK}
  \begin{align}
    &  \underline{q}^i= (1-c_i)\underline{q}^n + \sum_{j=1}^s a_{ij}\left(\underline{q}^n + \Delta t \underline{R}^j \right), \quad \mbox{for} \quad i= 1, ..., s,\\
    & \underline{q}^{RK} = \sum_{j=1}^s b_j \left( \underline{q}^n + \Delta t \underline{R}^j \right).
  \end{align}
\end{subequations} 
With this representation, the terms $\underline{q}^n + \Delta t \underline{R}^j$ explicitly appear within the RK formulation. If these terms, provided an SSP spatial discretization, are stable in a convex function $G$ for nonzero time steps, and the RK coefficients satisfy $0 \leq a_{ij}, c_i, b_j \leq 1$, then we will show that the RK solution and all intermediate stage solutions remain stable in $G$. Therefore, in this study we restrict our attention to RK methods with coefficients in the interval $[0,1]$. 

\begin{assumption}\label{Assumption_RK}
We assume that the RK methods considered in this work have coefficients satisfying $0 \leq a_{ij}, c_i, b_j \leq 1$.
\end{assumption}

\begin{remark}
Assumption \ref{Assumption_RK} covers non-SSP schemes such as the classical RK44 method.
\end{remark}

We now formally establish the connection between the stability of the terms $\underline{q}^n + \Delta t \underline{R}^j$ and the stability of RK methods. 

\begin{lemma}\label{lemma_SSP_modified_RK}
Assume that for a RK method satisfying Assumption \ref{Assumption_RK}, each stage derivative $\underline{R}^j$ satisfies $G \left( \underline{q}^n + \Delta t \underline{R}^j \right)\leq G\left( \underline{q}^n \right)$ with a nonzero time-step limit, where $G$ is a norm, semi-norm, or convex functional. Then, the RK solution and each stage solutions will be stable in $G$ with nonzero time-steps
\begin{equation}
G \left( \underline{q}^{j} \right) \leq G\left( \underline{q}^n \right), \quad \textit{and} \quad  G \left( \underline{q}^{RK} \right) \leq G\left( \underline{q}^n \right).
\end{equation}
\end{lemma}

\begin{proof}[Proof of Lemma \ref{lemma_SSP_modified_RK}]
For the considered RK methods, the coefficients satisfy $0 \leq 1-c_i, a_{ij}, b_j \leq 1$. Therefore, according to the RK representation \ref{Modified_RK}, the RK solution and all intermediate stage solutions can be expressed as a convex combination of $\underline{q}^n + \Delta t \underline{R}^j$ terms. With the convexity of $G\left(\underline{q} \right)$, the proof will be completed. 
\end{proof}

In the following subsections, we study the stability of the terms $\underline{q}^n + \Delta t \underline{R}^j$ for small time-steps. First, we demonstrate that these terms maintain stability with respect to strictly convex functionals, such as mathematical entropy. Then, we show that, when using the Lax-Friedrichs spatial discretization, these terms also preserve the positivity of density and internal energy. Finally, we demonstrate that, with the first-order upwind scheme and the second-order MUSCL spatial discretization, the terms $\underline{q}^n + \Delta t \underline{R}^j$ maintain the total variation diminishing property for Burgers' equation.

\subsection{Strictly convex functionals}

In this subsection, we assume that $G\left( \underline{q} \right)$ is a smooth, strictly convex functional, such as numerical entropy. Consequently, for any nonzero $R\left( \underline{q}(t) \right)$, 
\begin{equation}
R\left( \underline{q}(t) \right)^T \nabla^2 G\left(\underline{q}(t) \right) R\left( \underline{q}(t) \right) >0.
\end{equation}
We then show that if the forward Euler solution is stable with respect to $G$ for $\Delta t \leq \Delta t_{FE}$, the terms $\underline{q}^n + \Delta t \underline{R}^j$ also remain stable in $G$ provided that the time step is sufficiently small. 

\begin{lemma}\label{lemma_convex_functional}
Assume that $G$ is a strictly convex functional, and the forward Euler method is entropy stable for $\Delta t \leq \Delta t_{FE}$. Then, RK stage derivatives $\underline{R}^j$ satisfy
\begin{equation}
G\left(\underline{q}^n + \Delta t \underline{R}^j \right) \leq G\left(\underline{q}^n \right),
\end{equation}
when the step size is sufficiently small.
\end{lemma}

\begin{proof}[Proof of Lemma \ref{lemma_convex_functional}]
First assume that $\underline{R}^n= R\left( \underline{q}^n \right) \neq \underline{0}$. With the forward Euler solution, we have
\begin{equation}
G\left(\underline{q}^n + \Delta t \underline{R}^n \right)= G\left(\underline{q}^n \right) + \Delta t \nabla G\left(\underline{q}^n \right)^T\underline{R}^n + 0.5 \Delta t^2 {\underline{R}^n}^T \nabla^2 G\left(\underline{q}^n \right)\underline{R}^n  + \mathcal{O}(\Delta t^3).
\end{equation}
Since $ {\underline{R}^n}^T \nabla^2 G\left(\underline{q}^n \right)\underline{R}^n >0$, and the forward Euler solution is stable for small step sizes, it follows that $\nabla G\left(\underline{q}^n \right)^T\underline{R}^n <0$. On the other hand, we have
\begin{equation}
G \left( \underline{q}^n + \Delta t \underline{R}^j \right)= G\left(\underline{q}^n \right) + \Delta t \nabla G\left(\underline{q}^n \right)^T\underline{R}^n + \mathcal{O}(\Delta t^2).
\end{equation}
Since $\nabla G\left(\underline{q}^n \right)^T\underline{R}^n <0$, we will have $G\left(\underline{q}^n + \Delta t R^j \right) \leq G\left(\underline{q}^n \right)$ for small nonzero step sizes.

On the other hand, if $\underline{R}^n = \underline{0}$, the solution vector will not change and for each $j$ we will have $\underline{q}^j= \underline{q}^n$ and $\underline{R}^j= R \left( \underline{q}^n \right) = \underline{0}$, and consequently $G \left( \underline{q}^n + \Delta t \underline{R}^j \right)= G\left(\underline{q}^n \right)$.

\end{proof}

Then, the stability of RK methods satisfying the Assumption \ref{Assumption_RK} when step size is small enough can be concluded.

\begin{theorem}\label{Theorem_entropy_stability}
When $G$ is a strictly convex functional and the forward Euler method is stable in this norm for a range of positive step sizes $\Delta t \leq \Delta t_{FE}$, then, with RK methods satisfying Assumption \ref{Assumption_RK}, the RK solution and all stage solutions remain stable in $G$
\begin{equation}
G\left(\underline{q}^j \right) \leq G\left(\underline{q}^n \right), \quad  G\left(\underline{q}^{RK} \right) \leq G\left(\underline{q}^n \right).
\end{equation} 
\end{theorem}

\begin{proof}[Proof of Theorem \ref{Theorem_entropy_stability}]
The result directly follows from Lemmas \ref{lemma_SSP_modified_RK} and \ref{lemma_convex_functional}.
\end{proof}

\subsection{Positivity of density}

For the one-dimensional Euler equations, the vector of conserved variables at each point in space is given by
\begin{equation}
\underline{q}_i = \left[ \begin{array}{lcr}
\rho_i \\
m_i\\
E_i
\end{array} \right],
\end{equation}
with
\begin{equation}
m= \rho u, \quad \rho e=  E-  \frac{m^2}{2\rho}, \quad p= (\gamma -1 )\rho e,
\end{equation}
where $\rho$ is the density, $m$ is the momentum, $u$ is the velocity, $E$ is the total energy, $e$ is the internal energy, $p$ is the pressure, and $\gamma$ is the ratio specific heats. Admissible solutions for Euler equations must maintain positive density and internal energy. 

Positivity of density and pressure can be ensured, for example, by using the Lax-Friedrichs spatial discretization \cite{zhangPositivitypreservingHighOrder2010}
\begin{equation}
\frac{d}{dt} \underline{q}_i (t)= - \frac{1}{\Delta x} \left( h\left(\underline{q}_i , \underline{q}_{i+1} \right) - h\left(\underline{q}_{i-1} , \underline{q}_{i} \right)  \right),
\end{equation}
where $h$ is the Lax-Friedrichs flux
\begin{equation}
h\left( \underline{q}_r, \underline{q}_k \right)= \frac{1}{2} \left( f\left(\underline{q}_r \right) + f \left( \underline{q}_k \right) - a \left( \underline{q}_k - \underline{q}_r  \right) \right), \quad a= \max \left( |u_r| + \sqrt{\gamma \frac{p_r}{\rho_r}}, |u_k| + \sqrt{\gamma \frac{p_k}{\rho_k}} \right),
\end{equation}
and $f$ is the Euler flux function. When this semi-discrete system is advanced in time with the forward Euler method
\begin{equation}
\underline{q}_i^{FE}= (1 - \lambda a^n) \underline{q}_i^n + \frac{\lambda a^n}{2} \left[ \underline{q}_{i+1}^n - \frac{1}{a^n} f \left( \underline{q}_{i+1}^n \right)  \right] + \frac{\lambda a^n}{2} \left[ \underline{q}_{i-1}^n + \frac{1}{a^n} f \left( \underline{q}_{i-1}^n \right)  \right],
\end{equation}
with $\lambda= \Delta t/\Delta x$, density and internal energy remain positive provided that $\lambda a^n \leq 1$ \cite{zhangPositivitypreservingHighOrder2010}. We show in this part that with the Lax-Friedrichs method integrated with RK methods fulfilling Assumption \ref{Assumption_RK}, positivity of density can be maintained with nonzero time steps. As usual, we start by analyzing the terms $\underline{q}^n + \Delta t \underline{R}^j$.

\begin{lemma}\label{lemma_positivity_density}
With the Lax-Friedrichs spatial discretization, the density with $\underline{q}^n + \Delta t \underline{R}^j$ remains positive for small enough, yet nonzero, step sizes.

\end{lemma}

\begin{proof}[Proof of Lemma \ref{lemma_positivity_density}]
For the $\underline{q}^n + \Delta t \underline{R}^j$, we have
\begin{equation}
\underline{q}_i^n + \Delta t \underline{R}_i^j=  \underline{q}_i^n - \lambda a^j \underline{q}_i^j + \frac{\lambda a^j}{2} \left[ \underline{q}_{i+1}^j - \frac{1}{a^j} f \left( \underline{q}_{i+1}^j \right)  \right] + \frac{\lambda a^j}{2} \left[ \underline{q}_{i-1}^j + \frac{1}{a^j} f \left( \underline{q}_{i-1}^j \right)  \right],
\end{equation}
Since $\lambda a^j>0$ and the density components in the terms $\underline{q}_{i+1}^j - \frac{1}{a^j} f \left( \underline{q}_{i+1}^j \right)$ and $\underline{q}_{i-1}^j + \frac{1}{a^j} f \left( \underline{q}_{i-1}^j \right)$  are positive according to \cite{zhangPositivitypreservingHighOrder2010}, it remains to show that for small enough step sizes $\rho_i^n - \lambda a^j \rho_i^j >0$, where $\rho_i^n>0$. As $a^j= a^n + \mathcal{O}(\Delta t)$ and $\rho_i^j = \rho_i^n + \mathcal{O}(\Delta t)$
\begin{equation}
\rho_i^n - \lambda a^j \rho_i^j= \rho_i^n - \frac{\Delta t}{\Delta x} (a^n + \mathcal{O}(\Delta t))(\rho_i^n + \mathcal{O}(\Delta t))= \rho_i^n - \frac{\Delta t}{\Delta x} a^n \rho_i^n + \mathcal{O}(\Delta t^2).
\end{equation}
Thus, the condition $\rho_i^n - \lambda a^j \rho_i^j >0$ reduces to have
\begin{equation}
\frac{1}{a^n} > \mathcal{O}(\Delta t),
\end{equation}
which will be satisfied with small, nonzero step sizes. 
\end{proof}

\begin{theorem}\label{Theorem_positivity_density}
With the Lax-Friedrichs spatial discretization, RK methods satisfying the Assumption \ref{Assumption_RK} maintain positivity of the density with sufficiently small time steps.
\end{theorem}

\begin{proof}[Proof of Theorem \ref{Theorem_positivity_density}]
The proof follows directly from Lemmas \ref{lemma_SSP_modified_RK} and \ref{lemma_positivity_density}.  
\end{proof}

\begin{remark}
Lemma \ref{lemma_positivity_density} and Theorem \ref{Theorem_positivity_density} can be readily extended to the local Lax-Friedrichs spatial discretization.
\end{remark}

\subsection{Positivity of internal energy}

Assume that, for the Euler equations, a spatial discretization (e.g., Lax-Friedrichs) integrated with the forward Euler method preserves positivity of density and internal energy, and that the density with $\underline{q}^n + \Delta t \underline{R}^j$ terms remains positive with sufficiently small step sizes. We show here that the internal energy also remains positive with $\underline{q}^n + \Delta t \underline{R}^j$ terms with small enough step sizes.

\begin{lemma}\label{lemma_positivity_energy}
Assume that, with a given spatial discretization, the forward Euler method maintains positivity of density and internal energy for step sizes $\Delta t \leq \Delta t_{FE}$, with $\Delta t_{FE}>0$, and that the density with the terms $\underline{q}^n + \Delta t \underline{R}^j$ of an RK method remains positive with small step sizes. Then, internal energy also remains positive with $\underline{q}^n + \Delta t \underline{R}^j$ terms for sufficiently small time steps. 

\end{lemma}

\begin{proof}[Proof of Lemma \ref{lemma_positivity_energy}]
To start, for a node $i$ we define the following
\begin{equation}
\underline{q}_i^n = \begin{bmatrix}
\rho^n\\
m^n\\
E^n
\end{bmatrix}, \quad \underline{R}_i^n= \begin{bmatrix}
R^{\rho n}\\
R^{mn}\\
R^{En}
\end{bmatrix}, \quad \underline{R}_i^j= \begin{bmatrix}
R^{\rho j}\\
R^{mj}\\
R^{Ej}
\end{bmatrix}.
\end{equation}
The internal energy at point $i$ with the forward Euler solution is
\begin{equation}
\rho e \left( \underline{q}^{FE} \right)=  \frac{(E^n + \Delta t R^{En})(\rho^n + \Delta t R^{\rho n}) - 0.5(m^n + \Delta t R^{mn})^2}{\rho^n + \Delta t R^{\rho n}},             
\end{equation}
where both the numerator and denominator remain positive for step sizes satisfying $0 \leq \Delta t \leq \Delta t_{FE}$. In particular, the numerator is a quadratic polynomial in $\Delta t$
\begin{equation}\label{quadratic_eq_energy}
 c_{n,n} + b_{n,n} \Delta t + a_{n,n} \Delta t^2 ,                          
\end{equation}
where
\begin{equation}
a_{n,n}= R^{En} R^{\rho n} - 0.5\left(R^{mn} \right)^2, 
\end{equation}
\begin{equation}
b_{n,n}= E^n R^{\rho n} + \rho^n R^{En} - m^n R^{mn}
\end{equation}
\begin{equation}
c_{n,n}= E^n \rho^n - 0.5\left(m^n \right)^2. 
\end{equation}
We have $c_{n,n}= \rho\left( \underline{q}^{FE} \right) \rho e \left( \underline{q}^{FE} \right) >0$, and Eq. \ref{quadratic_eq_energy} has no roots in the interval $0 \leq \Delta t \leq \Delta t_{FE}$. 

On the other hand, for $\underline{q}^n + \Delta t \underline{R}^j$ we have
\begin{equation}\label{energy_eq_j}
\rho e \left(\underline{q}^n + \Delta t \underline{R}^j \right)= \frac{(E^n + \Delta t R^{Ej})(\rho^n + \Delta t R^{\rho j}) - 0.5(m^n + \Delta t R^{mj})^2}{\rho^n + \Delta t R^{\rho j}} ,
\end{equation}
where here the denominator by lemma's assumption remains positive for small step sizes, and the numerator is no longer a quadratic polynomial, as $\underline{R}^j= \underline{R}^n + \mathcal{O}(\Delta t)$. However, the the numerator of Eq. \ref{energy_eq_j} converges to the value of Eq. \ref{quadratic_eq_energy} as $\Delta t \rightarrow 0$. Since Eq. \ref{quadratic_eq_energy} on $[0, \Delta t_{FE}]$ is positive, there should exist a positive time step limit $\Delta t^*>0$ such that the numerator of Eq. \ref{energy_eq_j} remains positive for $\Delta t \in [0, \Delta t^*]$.

\end{proof}

We can now extend the positivity of internal energy with $\underline{q}^n + \Delta t \underline{R}^j$ terms to the RK solution and intermediate stage solutions of RK methods satisfying the Assumption \ref{Assumption_RK}.

\begin{theorem}\label{Theorem_positivity_energy}

When forward Euler time integration maintains positivity of density and internal energy for a given spatial discretization, and the density with the terms $\underline{q}^n + \Delta t \underline{R}^j$ of an RK method remains positive for small step sizes, then the internal energy also remains positive for sufficiently small step sizes.
\end{theorem}

\begin{proof}[Proof of Theorem \ref{Theorem_positivity_energy}]
By using the Lemmas \ref{lemma_SSP_modified_RK} and \ref{lemma_positivity_energy}, the proof can be completed. 
\end{proof}

\subsection{Total variation}

For scalar hyperbolic conservation laws, such as Burgers' equation, the Total Variation of a discrete solution is defined as 
\begin{equation}
TV \left(\underline{q} \right)= \sum_i |q_{i+1} - q_i|,
\end{equation}
and a spatial discretization is called Total Variation Diminishing (TVD) if, when integrated using the forward Euler method,
\begin{equation}
TV \left(\underline{q}^n + \Delta t \underline{R}^n \right) \leq TV\left(\underline{q}^n \right),
\end{equation}
for step sizes $\Delta t \leq \Delta t_{FE}$, where $\Delta t_{FE}>0$. In this section, we demonstrate that RK methods satisfying Assumption \ref{Assumption_RK} can maintain the TVD property of certain spatial discretizations with nonzero step sizes.

To proceed, we introduce a lemma concerning the TV of the terms $\underline{q}^n + \Delta t \underline{R}^j$ in the special case where $q_i^n \neq q_{i+1}^n$ for all $i$. This lemma will be used in the proofs presented in the following subsections.

\begin{lemma}\label{Lemma_TVD_first}
Let $\underline{q}^n$ be a solution vector such that $q_i^n \neq q_{i+1}^n$ for all $i$. Suppose that for every local minimum degree of freedom $\ell$, i.e., $q_{\ell}^n < \{ q_{\ell-1}^n,\; q_{\ell+1}^n \}$, we have $R_{\ell}^j\geq 0$, and for every local maximum degree of freedom $m$, i.e., $ q_m^n > \{q_{m-1}^n, \; q_{m+1}^n \}$, we have $R_m^j \leq 0$ for sufficiently small step sizes. Then $TV\left(\underline{q}^n + \Delta t \underline{R}^j \right) \leq TV \left(\underline{q}^n \right)$.
\end{lemma}

\begin{proof}[Proof of Lemma \ref{Lemma_TVD_first}]
Assume that for the solution vector $\underline{q}^n$, indices $\ell$, $m$, and $k$ represent consecutive local extrema, with no other local extrema between them
\begin{equation}
\underline{q}^n= \left[\dots, q_{\ell}^n, \dots, q_m^n, \dots, q_k^n, \dots \right]^T,
\end{equation}
where among each pair of consecutive local extrema, one is a local maximum and the other one is a local minimum. The TV of $\underline{q}^n$ then can be written as
\begin{equation}
TV \left(\underline{q} \right)= \dots + |q_m^n - q_{\ell}^n | \;+ \; |q_k^n - q_m^n| + \dots \, .
\end{equation}
Therefore, for the vector $\underline{q}^n$, the TV is the sum of absolute differences between adjacent local extrema.

Without loss of generality, assume that $\ell$ is a local minimum. Then, $m$ must be a local maximum, and $k$ a local minimum. With $\underline{q}^n +\Delta t \underline{R}^j$, provided that the time step size remains small, any local minimum degree of freedom $\ell$ remains a local minimum, since there exist $\Delta t>0$ such that
\begin{equation}
q_{\ell}^n < \left\{q_{\ell-1}^n,\; q_{\ell+1}^n \right\}, \quad q_{\ell}^n + \Delta t R_{\ell}^j < \left\{q_{{\ell}-1}^n+ \Delta t R_{{\ell}-1}^j,\; q_{{\ell}+1}^n+ \Delta t R_{{\ell}+1}^j \right\}.
\end{equation}
By a similar argument, for sufficiently small step sizes, $m$ remains a local maximum, $k$ remains a local minimum, and no new local extrema appear among the degrees of freedom $i \in (\ell,m)$ and $i \in (m,k)$. Therefore
\begin{equation}
TV \left(\underline{q}^n+ \Delta t \underline{R}^j \right)= \dots + |q_m^n + \Delta t R_m^j - q_{\ell}^n - \Delta t R_{\ell}^j | \;+ \; |q_k^n + \Delta t R_k^j - q_m^n -\Delta t R_m^j | + \dots \, .
\end{equation}

Since $\ell, m$, and $k$ are consecutive local extrema, for small step sizes 
\begin{equation}
q_m^n + \Delta t R_m^j > \left\{ q_{\ell}^n + \Delta t R_{\ell}^j, \; q_k^n + \Delta t R_k^j  \right\}.
\end{equation}
Moreover, by assumption, we have $R_{\ell}^j\geq 0$, $R_m^j\leq 0$, and $R_k^j\geq 0$. Therefore
\begin{equation}
|q_m^n + \Delta t R_m^j - q_{\ell}^n - \Delta t R_{\ell}^j | \leq |q_m^n  - q_{\ell}^n  |, \quad |q_k^n + \Delta t R_k^j - q_m^n - \Delta t R_m^j |\leq |q_k^n  - q_m^n  |.
\end{equation}
Consequently
\begin{equation}
TV\left(\underline{q}^n + \Delta t \underline{R}^j \right) \leq TV \left(\underline{q}^n \right).
\end{equation}
\end{proof}

\subsubsection{First-order upwind method}
With the first-order upwind method, Burgers' equation is discretized as \cite{levequeFiniteVolumeMethods2002}
\begin{equation}\label{first_order_upwind}
 \frac{d}{dt}q_i(t)= - \frac{1}{\Delta x} \left( f(q_i ) - f (q_{i-1} ) \right)             , \quad f\left(q_i \right)= \frac{q_i^2}{2}.
\end{equation}
According to Harten's lemma \cite{hartenHighResolutionSchemes1982}, when the solution $q_i(t)$ remains within the range $[0,1]$, this spatial discretization is TVD with $\Delta t_{FE}= \Delta x$. In the following, we demonstrate the TVD property of $\underline{q}^n + \Delta t \underline{R}^j$ terms when RK methods satisfying Assumption \ref{Assumption_RK} are employed.

\begin{lemma}\label{lemma_TVD_first_order_upwind}
When the ODE system \ref{first_order_upwind} with $0 \leq q_i^n \leq 1$ is integrated using RK methods satisfying Assumption \ref{Assumption_RK}, we will have $TV\left(\underline{q}^n + \Delta t \underline{R}^j \right) \leq TV \left(\underline{q}^n \right)$ for $j= 1, \dots, s$, when the step size is sufficiently small, yet nonzero.
\end{lemma}

\begin{proof}[Proof of Lemma \ref{lemma_TVD_first_order_upwind}]

First, assume that $q_i^n \neq q_{i+1}^n$ for any $i$. For any local minimum degree of freedom $\ell$, it is easy to show that $R_{\ell}^n>0$ and
\begin{equation}
R_{\ell}^j= R_{\ell}^n + \mathcal{O}(\Delta t)>0,
\end{equation}
for small step sizes. Similarly, for any local maximum degree of freedom $m$ we have $R_m^n<0$, and 
\begin{equation}
 R_m^j= R_m^n + \mathcal{O}(\Delta t)<0,
\end{equation}
again for small step sizes. Therefore, according to Lemma \ref{Lemma_TVD_first}, $TV \left(\underline{q}^n + \Delta t \underline{R}^j \right)\leq TV \left(\underline{q}^n \right)$.

Next, assume that for some $i$, $q_{i-1}^n> q_{i}^n= q_{i+1}^n$. Under this condition, with the assumed RK methods and small nonzero step sizes, it can be shown that
\begin{equation}
q_{i-1}^n + \Delta t R_{i-1}^j > q_i^n + \Delta t R_i^j \geq q_{i+1}^n + \Delta t R_{i+1}^j \geq q_{i+1}^n.
\end{equation}
More generally, when $q_{i-1}^n > q_i^n = ...= q_{i+r}^n$ with $r\geq 1$, we will have 
\begin{equation}
q_{i-1}^n + \Delta t R_{i-1}^j> q_i^n + \Delta t R_i^j \geq ... \geq q_{i+r}^n + \Delta t R_{i+r}^j \geq q_{i+r}^n.
\end{equation}
Therefore, degrees of freedom $i$ to $i+r-1$ do not contribute to TV of $\underline{q}^n + \Delta t \underline{R}^j$, .i.e., 
\begin{equation}
TV\left( \underline{q}^n + \Delta t \underline{R}^j\right)= TV([..., q_{i-1}^n + \Delta t R_{i-1}^j, q_{i+r}^n + \Delta t R_{i+r}^j, ...]^T).
\end{equation}

Similarly, when for some $i$, $q_{i-1}^n < q_i^n = ...= q_{i+r}^n$ with $r \geq 1$, with the considered RK methods and small step sizes
\begin{equation}
q_{i-1}^n + \Delta t R_{i-1}^j< q_i^n + \Delta t R_i^j \leq ... \leq q_{i+r}^n + \Delta t R_{i+r}^j \leq q_{i+r}^n.
\end{equation}
Once again, the degrees of freedom $i$ to $i+r-1$ do not contribute to the TV of $\underline{q}^n + \Delta t \underline{R}^j$.

As a result, when $q_{i-1}^n \neq q_i^n= ...=q_{i+r}^n \neq q_{i+r+1}^n$ with $r\geq 1$, we can consider a modified solution vector for calculating the TV
\begin{equation}
TV\left(\underline{q}^n \right)= TV \left( \left[\dots, q_{i-1}^n, q_{i+r}^n, \dots \right]^T \right),
\end{equation}
and
\begin{equation}
\begin{aligned}
TV\left(\underline{q}^n + \Delta t \underline{R}^j\right) &= TV \left(\left[\dots, q_{i-1}^n + \Delta t R_{i-1}^j, q_{i+r}^n + \Delta t R_{i+r}^j , \dots \right]^T \right)\\
&= TV \left( \left[\dots, q_{i-1}^n, q_{i+r}^n,\dots \right]^T + \Delta t \left[\dots, R_{i-1}^j, R_{i+r}^j, \dots \right]^T \right),
\end{aligned}
\end{equation}
where, for the vector $[..., q_{i-1}^n, q_{i+r}^n,...]^T$, no consecutive degrees of freedom are equal. With this modified vector, every local minimum degree of freedom $\ell$ satisfies $R_{\ell}^j \geq 0$, and every local maximum degree of freedom $m$ satisfies $R_m^j \leq 0$. Therefore, according to Lemma \ref{Lemma_TVD_first}
\begin{equation}
TV\left(\underline{q}^n + \Delta t \underline{R}^j \right)= TV \left( \left[\dots, q_{i-1}^n + \Delta t R_{i-1}^j, q_{i+r}^n + \Delta t R_{i+r}^j,\dots \right]^T \right) \leq TV\left( \left[\dots, q_{i-1}^n,q_{i+r}^n, \dots \right]^T \right)= TV\left(\underline{q}^n \right).
\end{equation}

\end{proof}

Then, we can show the TVD property of RK methods applied to Burgers' equation discretized with the first-order upwind scheme.

\begin{theorem}\label{Theorem_TVD_upwind}
When Burgers' equation is discretized with the first-order upwind scheme and $0\leq q_i^n \leq 1$, time integration with RK methods satisfying Assumption \ref{Assumption_RK} results in maintaining the TVD property by RK stage solutions and the step solution for nonzero small steps
\begin{equation}
TV\left(\underline{q}^j \right) \leq TV\left(\underline{q}^n\right), \quad \textit{and} \quad TV\left(\underline{q}^{RK} \right) \leq TV\left(\underline{q}^n\right).
\end{equation}
\end{theorem}

\begin{proof}[Proof of Theorem \ref{Theorem_TVD_upwind}]
By using Lemmas \ref{lemma_SSP_modified_RK} and \ref{lemma_TVD_first_order_upwind}, the result will be obtained. 
\end{proof}

\subsubsection{Second-order MUSCL scheme}

Next, we consider the standard minmod based MUSCL second-order spatial discretization of Burgers' equation \cite{gottliebTotalVariationDiminishing1998}
\begin{equation}
\frac{d}{dt}q_i (t) = - \frac{1}{\Delta x} \left(  \hat{f}_{i+ \frac{1}{2}} - \hat{f}_{i-\frac{1}{2}} \right),
\end{equation} 
where \begin{equation}
\hat{f}_{i+ \frac{1}{2}} = h \left(  q_{i+ \frac{1}{2}}^{-} ,q_{i+ \frac{1}{2}}^{+}  \right),
\end{equation}
with
\begin{equation}
 q_{i+ \frac{1}{2}}^{-} = q_i + \frac{1}{2} \textit{minmod}\left( q_{i+1} - q_i, q_i - q_{i-1}  \right),
\end{equation}
\begin{equation}
 q_{i+ \frac{1}{2}}^{+} = q_{i+1} - \frac{1}{2} \textit{minmod}\left( q_{i+2} - q_{i+1}, q_{i+1} - q_{i}  \right),
\end{equation}
where $h$ is the Gudunov flux 

\begin{equation}
h\left( q^-, q^+   \right)= \left\{  \begin{array}{lcr}
\min_{q^- \leq q \leq q^+} \left( \frac{q^2}{2} \right)  & \textit{if} \; q^- \leq q^+ ,\\
\max_{q^- \geq q \geq q^+}\left( \frac{q^2}{2} \right) & \textit{if} \; q^- > q^+ ,
\end{array} \right.
\end{equation} 
and 
\begin{equation}
\textit{minmod}(a,b)= \frac{\textit{sign}(a) + \textit{sign}(b)}{2} \min(|a|, |b|).
\end{equation}
This spatial discretization is TVD with $\Delta t_{FE}= \frac{\Delta x}{2 \max_i |q_i^n|}$ \cite{gottliebTotalVariationDiminishing1998}. We show that by using RK methods satisfying Assumption \ref{Assumption_RK}, $TV \left(\underline{q}^n + \Delta t \underline{R}^j \right) \leq TV \left(\underline{q}^n \right)$, $j=1, \dots, s$, when step size is small enough.

\begin{lemma}\label{lemma_TVD_minmod}
When the standard minmod based second-order MUSCL spatial discretization is integrated with RK methods satisfying the Assumption \ref{Assumption_RK}, we will have $TV\left(\underline{q}^n + \Delta t \underline{R}^j \right) \leq TV\left(\underline{q}^n\right)$ for small nonzero time steps. 
\end{lemma}

\begin{proof}[Proof of Lemma \ref{lemma_TVD_minmod}]

We begin by assuming that $q_i^n \neq q_{i+1}^n$ for all $i$. For any local minimum degree of freedom $\ell$, $q_{\ell}^j$ will also be a local minimum when $\Delta t$ is small enough, and to calculate $R_{\ell}^j$ we have
\begin{equation}
q_{\ell+1/2}^{j-}= q_{\ell}^j, \quad q_{\ell+1/2}^{j+} > q_{\ell}^j, \quad \hat{f}_{\ell+ \frac{1}{2}}^j \leq \frac{\left(q_{\ell}^j\right)^2}{2},
\end{equation}
and
\begin{equation}
q_{\ell-\frac{1}{2}}^{j-} > q_{\ell}^j, \quad q_{\ell-\frac{1}{2}}^{j+}= q_{\ell}^j, \quad \hat{f}_{\ell-\frac{1}{2}}^j \geq \frac{\left(q_{\ell}^j\right)^2}{2}.
\end{equation}
Therefore
\begin{equation}
R_{\ell}^j= - \frac{1}{\Delta x} \left( \hat{f}_{\ell+\frac{1}{2}}^j - \hat{f}_{\ell-\frac{1}{2}}^j \right) \geq 0.
\end{equation}
With a similar approach, it can be shown that for any locally maximum degree of freedom $m$, we have $R_m^j \leq 0$. Therefore, according to Lemma \ref{Lemma_TVD_first}, $TV \left(  \underline{q}^n + \Delta t \underline{R}^j \right) \leq TV \left( \underline{q}^n \right)$.

Then, when for some $i$, $q_{i-1}^n> q_i^n=\dots= q_{i+r}^n$, with considered RK methods and small nonzero step sizes, it can be shown that
\begin{equation}
q_{i-1}^n + \Delta t R_{i-1}^j > q_i^n + \Delta t R_i^j \geq \dots \geq q_{i+r}^n + \Delta t R_{i+r}^j \geq q_{i+r}^n.
\end{equation}
Moreover, when $q_{i-1}^n < q_i^n=\dots= q_{i+r}^n$, with small step sizes 
\begin{equation}
q_{i-1}^n + \Delta t R_{i-1}^j < q_i^n + \Delta t R_i^j \leq \dots \leq q_{i+r}^n + \Delta t R_{i+r}^j \leq q_{i+r}^n.
\end{equation}
Therefore, when $q_{i-1}^n \neq q_i^n=\dots = q_{i+r}^n \neq q_{i+r+1}^n$, degrees of freedom $i$ to $i+r-1$ do not contribute to the TV of $\underline{q}^n + \Delta t \underline{R}^j$
\begin{equation}
TV\left(\underline{q}^n + \Delta t \underline{R}^j \right)= TV \left( \left[ \dots, q_{i-1}^n + \Delta t R_{i-1}^j , q_{i+r}^n + \Delta t R_{i+r}^j, \dots \right] \right).
\end{equation}
With the solution vector $[..., q_{i-1}^n, q_{i+r}^n, ...]$, no consecutive solutions are equal, and for any local minimum degree of freedom $\ell$ we have $R_{\ell}^j \geq 0$ and for any local maximum $m$ we have $R_m^j \leq 0$. Therefore, according to Lemma \ref{Lemma_TVD_first}
\begin{equation}
TV \left(\underline{q}^n + \Delta t \underline{R}^j \right)= TV \left( \left[ \dots, q_{i-1}^n + \Delta t R_{i-1}^j, q_{i+r}^n + \Delta t R_{i+r}^j, \dots \right]^T  \right) \leq TV \left( \left[\dots, q_{i-1}^n,q_{i+r}^n, \dots  \right]^T \right)= TV\left(\underline{q}^n \right).
\end{equation}

\end{proof}

Now we can provide the main theorem of this part.

\begin{theorem}\label{Theorem_TVD_minmod}
When integrating the standard minmod based second-order MUSCL spatial discretization of Burgers' equation with RK methods falling under the Assumption \ref{Assumption_RK}, the TVD property will be maintained by the intermediate stage solutions and the RK step solution when the step size is small enough
\begin{equation}
TV\left(\underline{q}^j \right) \leq TV\left(\underline{q}^n\right), \quad \textit{and} \quad TV\left(\underline{q}^{RK} \right) \leq TV\left(\underline{q}^n\right).
\end{equation}
\end{theorem}

\begin{proof}[Proof of Theorem \ref{Theorem_TVD_minmod}]
Combining Lemmas \ref{lemma_SSP_modified_RK} and \ref{lemma_TVD_minmod} leads to the desired result.

\end{proof}

\section{Example cases}\label{Section_examples}

In this section, we test a series of RK methods that satisfy Assumption \ref{Assumption_RK}: the first-order forward Euler method (SSP), the second-order midpoint rule \cite{butcherNumericalMethodsOrdinary2008} (non-SSP), the third-order three-stage SSPRK33 scheme \cite{gottliebHighOrderStrong2009a} (SSP), the third-order three-stage RK31 method \cite{butcherNumericalMethodsOrdinary2008} (non-SSP), and the classical fourth-order RK44 method (non-SSP). These RK methods are applied to a series of challenging numerical problems that require preserving the stability of the forward Euler method. While the SSP time-step limit for non-SSP schemes and $\underline{q}^n + \Delta t \underline{R}^{j\neq 1}$ terms is exactly zero, we demonstrate that these schemes maintain stability in these example cases with nonzero step sizes. To quantify the practical time step limits, we report two parameters, $c^p$ and $c^s$, for each RK method in each test case. The value $c^p \Delta t_{FE}$ denotes the largest time step for which the RK solution and all intermediate stage solutions maintain the stability of the forward Euler method throughout the simulation, while $c^s \Delta t_{FE}$ is the largest time step for which the terms $\underline{q}^n + \Delta t \underline{R}^j$, $j\in \{1, \dots, s \}$, maintain forward Euler stability. We show in practice that for these experiments, $c^p$ and $c^s$ values for each of the RK methods are higher than zero, confirming that non-SSP RK methods and $\underline{q}^n + \Delta t \underline{R}^j$ terms maintained stability of the forward Euler method with nonzero step sizes, as expected from the theoretical findings.

\subsection{Energy dissipative Burgers' problem}
We consider Burgers' problem, with a periodic domain $x\in [-1, 1]$ with an initial condition 
\begin{equation}
q(x,0)= e^{-30x^2},
\end{equation}
discretized with an energy dissipative spatial discretization \cite{ketchesonRelaxationRungeKuttaMethods2019, najafianQuasiorthogonalRungeKuttaProjection2025}
\begin{equation}
\frac{d}{dt}q_i= -\frac{1}{\Delta x} (F_{i+1/2} - F_{i-1/2}), \quad F_{i+1/2}= \frac{q_i^2 + q_iq_{i+1} + q_{i+1}^2}{6} - \mu \left(q_{i+1} - {q}_i \right),             
\end{equation}
where $\mu = 1.E-3$. With this discretization, the energy $G\left( \underline{q} \right)= 0.5 \underline{q}^T \underline{q}$, which is a strictly convex function of the solution, decreases monotonically at each time step when the forward Euler method is used with $\Delta t \leq \Delta t_{FE}= 0.006 \Delta x$.

\begin{figure}[htbp]
\begin{subfigure}[b]{0.4\textwidth}
\centering
\includegraphics[width=\textwidth]{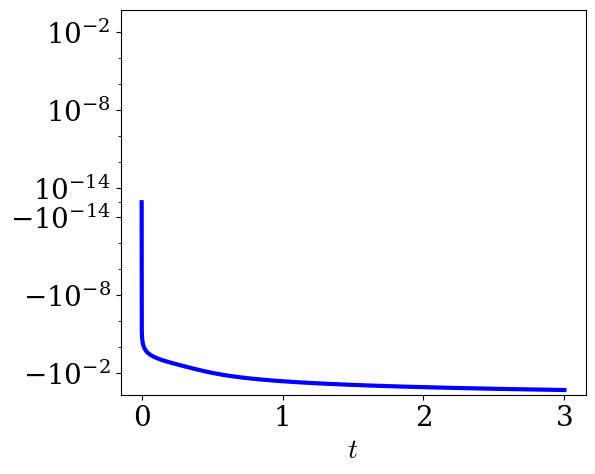}
\caption{$G\left(\underline{q}^{RK} \right) - G \left(\underline{q}^0\right)$}
\end{subfigure}
\begin{subfigure}[b]{0.4\textwidth}
\centering
\includegraphics[width=\textwidth]{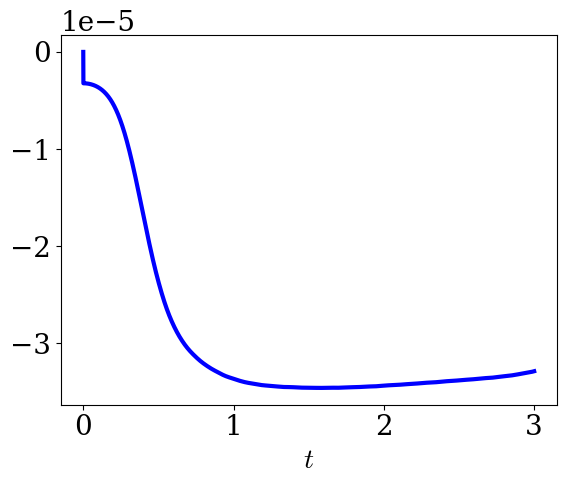}
\caption{$G\left(\underline{q}^{n} + \Delta t \underline{R}^{j=s} \right) - G\left( \underline{q}^n \right)$}
\end{subfigure}
\caption{Energy dissipative Burgers' problem integrated with RK44 using $\Delta t= \Delta t_{FE}$. This figure shows that both RK solution and $\underline{q}^n + \Delta t \underline{R}^j$ term preserved energy dissipative behavior of forward Euler method with a nonzero step size.}
\label{Fig_Burgers_dissipative_RK44}
\end{figure}

Figure \ref{Fig_Burgers_dissipative_RK44} demonstrates that, with a time step equal to $\Delta t_{FE}$, the RK44 solution and its $\underline{q}^n + \Delta t \underline{R}^{j=s}$ term preserveed the energy decay property of forward Euler method. Moreover, Table \ref{table_Burgers_dissipative} shows that $c^s$ and $c^p$ values for each of the tested RK methods are greater than zero, even for non-SSP schemes, confirming the theory developed in this work regarding stability in strictly convex functionals.

\begin{table}[htbp]
\centering
\begin{tabular}{m{3cm} m{1cm} m{1cm} m{1cm}}
\hline
RK scheme & $c^{ssp}$ & $c^s$ & $c^p$ \\ [0.5ex] 
\hline
forward Euler & $1$ &$1.0$  & $1.0$  \\
midpoint rule & $0$ & $1.0$ & $2.0$  \\
SSPRK33 & $1$ & $1.0$ & $1.0$  \\
RK31 & $0$ & $1.0$ & $1.5$ \\
RK44 & $0$ & $1.0$ & $2.0$ \\
\hline
\end{tabular}
\caption{$c^s$ and $c^p$ values for the dissipative Burgers' problem.}
\label{table_Burgers_dissipative}
\end{table}

\subsection{Burgers' problem with a first-order upwind method}

Consider the test case presented by Gottlieb et al. \cite{gottliebHighOrderStrong2009a}, involving Burgers' equation on the domain $x\in [0,2]$ with the initial condition
\begin{equation}
q(x,0)= \frac{1}{2} - \frac{1}{4}\sin(\pi x),
\end{equation}
and the periodic boundary condition. The solution will be right-traveling, and a shock will be developed over time. Using a first-order upwind finite difference spatial discretization, this problem will be TVD with the forward Euler method for $\Delta t \leq \Delta x$ \cite{gottliebHighOrderStrong2009a}, and the solution remains smooth even after the shock formation. This test case was used to demonstrate the advantage of using SSP time integration methods.

\begin{figure}[htbp]
\begin{subfigure}[b]{0.4\textwidth}
\centering
\includegraphics[width=\textwidth]{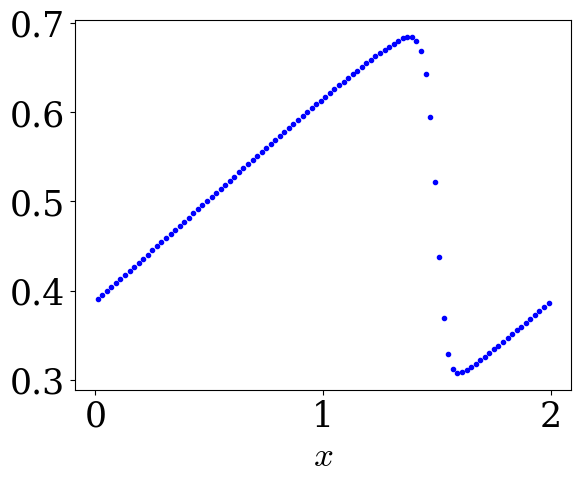}
\caption{Solution}
\end{subfigure}
\begin{subfigure}[b]{0.4\textwidth}
\centering
\includegraphics[width=\textwidth]{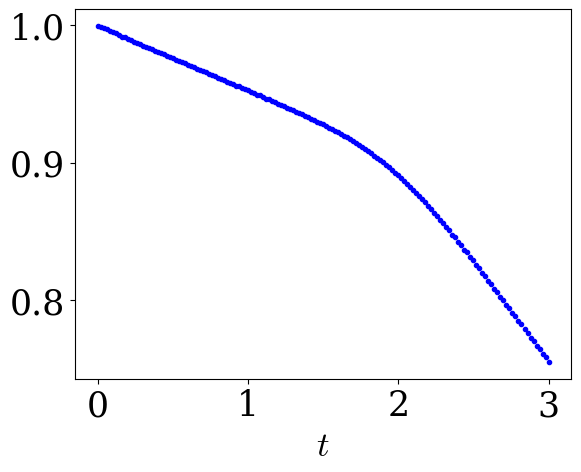}
\caption{$TV\left(\underline{q}^{RK} \right)$}
\end{subfigure}
\caption{Burgers' problem with a first-order upwind spatial discretization, integrated with RK44 using $\Delta t= \Delta t_{FE}$. It demonstrates that the TVD property is maintained with this non-SSP method using a nonzero time-step.}
\label{FD_RK44}
\end{figure}

\begin{figure}[htbp]
\centering
\includegraphics[width=0.5\textwidth]{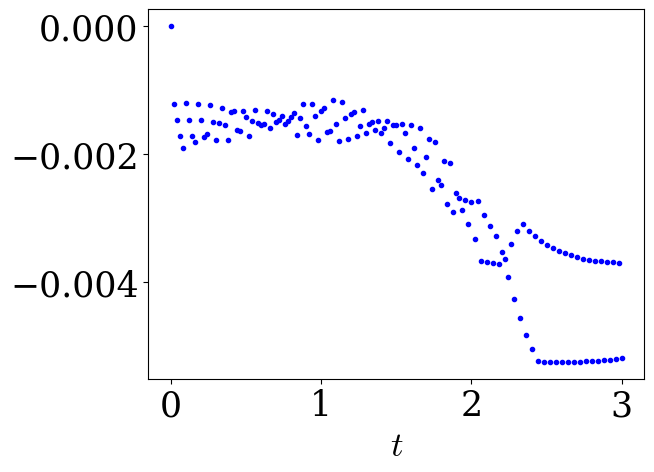}
\caption{ $TV\left(\underline{q}^n + \Delta t \underline{R}^{j=s} \right) - TV \left(\underline{q}^n \right)$ values of the RK44 scheme with $\Delta t= \Delta t_{FE}$, for the Burgers' problem with a first-order upwind spatial discretization. It shows that $\underline{q}^n + \Delta t \underline{R}^{j=s}$ terms maintain the TVD property of forward Euler method with a nonzero time-step . }
\label{FD_TVD_RK44_s}
\end{figure}

According to the Theorem \ref{Theorem_TVD_upwind}, RK methods under Assumption \ref{Assumption_RK} are expected to maintain TVD stability of forward Euler method in this problem. To verify this, the problem is integrated using RK44 with $\Delta t= \Delta t_{FE}$ and $\Delta x= 0.02$ up to $T=3$, which is beyond the point of shock formation. Figure \ref{FD_RK44} presents the RK44 solution at $T=3$, and the TV evolution  during the simulation, confirming that the solution remains smooth and TV monotonically decreases even with a step size equal to $\Delta t_{FE}$. Moreover, Figure \ref{FD_TVD_RK44_s} shows that $TV\left(\underline{q}^n + \Delta t \underline{R}^{j=s} \right) - TV \left(\underline{q}^n \right)$ remains negative throughout the simulation, confirming that $\underline{q}^n + \Delta t \underline{R}^{j=s}$ term maintains the TVD property, as predicted from the theory. In addition, Table \ref{table_Burgers_TVD} provides, for each RK method in this test case, $c^s$ and $c^p$ values. These values, which are greater than zero, confirm that the TVD property is maintained with considered non-SSP RK methods and $\underline{q}^n + \Delta t \underline{R}^j$ terms with positive time-steps.

\begin{table}[htbp]
\centering
\begin{tabular}{m{3cm}  m{1cm} m{1cm} m{1cm}}
\hline
RK scheme & $c^{ssp}$ &   $c^s$ & $c^p$ \\ [0.5ex] 
\hline
forward Euler & $1$   &$1.3$  & $1.3$  \\
midpoint rule & $0$   & $1.3$ & $1.6$  \\
SSPRK33 & $1$ &  $1.3$ & $1.3$  \\
RK31 & $0$    & $1.3$ & $2.0$ \\
RK44 & $0$ &  $1.3$ & $2.2$ \\
\hline
\end{tabular}
\caption{TVD time step limits for the Burgers' problem with a first-order upwind spatial discretization.}
\label{table_Burgers_TVD}
\end{table}

\subsection{Burgers' problem with a second order MUSCL method}

In this example, we follow the test case provided in \cite{gottliebTotalVariationDiminishing1998} which was used to show it is safer to use an SSP time integration scheme. This example considers Burgers' equation on the domain $x\in [-10, 70]$ with the initial condition 
\begin{equation}
q(x,0)= \left\{ \begin{array}{lcr}
1, \quad &\mbox{if} \; x\leq 0, \\
-0.5, \quad  &\mbox{if} \; x>0, \end{array} \right.
\end{equation}
discretized using the standard minmod-based MUSCL second-order spatial discretization with a uniform grade spacing $\Delta x=1$. With this discretization, the forward Euler method is TVD with a step size limit of $\Delta t_{FE}= \frac{\Delta x}{2 \max_i |q_i|}$. In \cite{gottliebTotalVariationDiminishing1998}, a non-SSP scheme with negative coefficients, thus violating Assumption \ref{Assumption_RK}, is tested to show the violation of the TVD stability by a non-SSP schemes. 

\begin{figure}[htbp]
\begin{subfigure}[b]{0.4\textwidth}
\centering
\includegraphics[width=\textwidth]{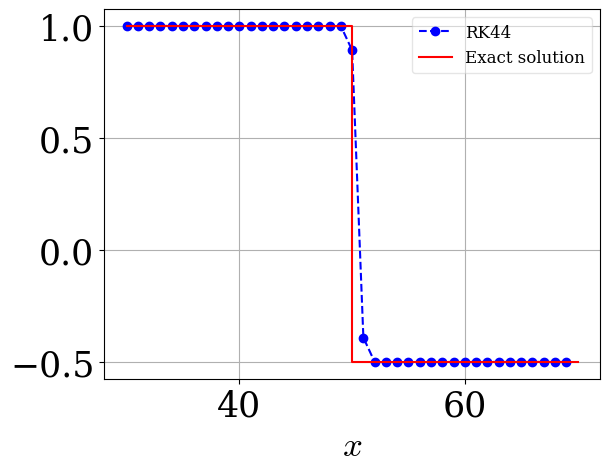}
\caption{Solution}
\end{subfigure}
\begin{subfigure}[b]{0.4\textwidth}
\centering
\includegraphics[width=\textwidth]{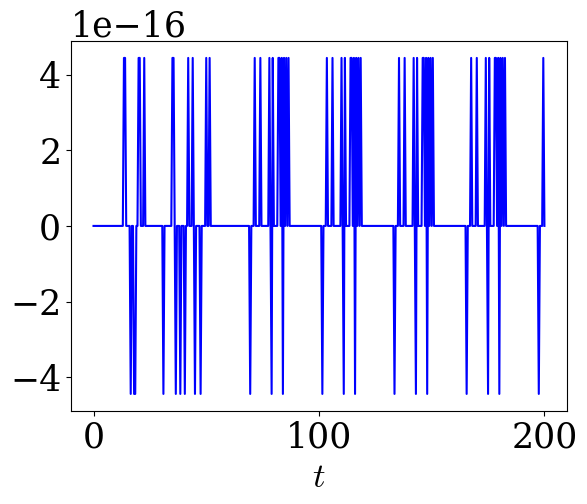}
\caption{$TV\left(\underline{q}^{RK}  \right) - TV\left( \underline{q}^0 \right)$}
\end{subfigure}
\caption{Burgers' problem discretized with the second-order MUSCL scheme, integrated with RK44 using $\Delta t= \Delta t_{FE}$, showing that RK44 maintains the TVD property of the forward Euler solution.}
\label{MUSCL2_RK44}
\end{figure}

\begin{figure}[htbp]
\centering
\includegraphics[width=0.4\textwidth]{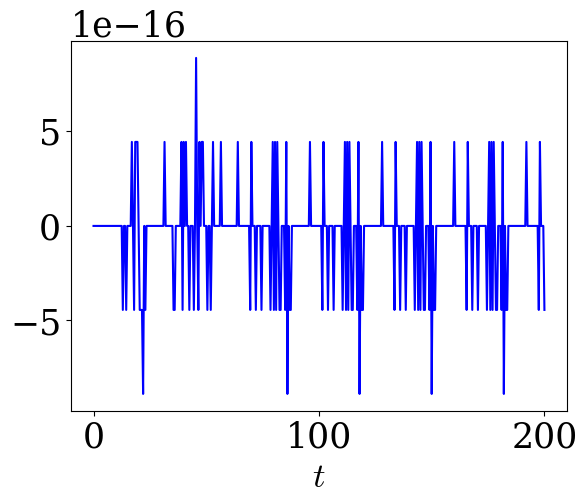}
\caption{$TV\left(\underline{q}^n + \Delta t \underline{R}^{j=s} \right) - TV \left(\underline{q}^n \right)$, when the Burgers' problem is discretized with the second-order MUSCL scheme and integrated with RK44 using $\Delta t= \Delta t_{FE}$.}
\label{MUSCL2_TVD_RK44_s}
\end{figure}

However, we have proved that RK44 and other non-SSP schemes satisfying the Assumption \ref{Assumption_RK} maintain the TVD property for nonzero step sizes. Figure \ref{MUSCL2_RK44} demonstrates that integrating the problem with RK44, using $\Delta t= \Delta t_{FE}$ up to $T= 200$, maintained the TVD property up to machine precision. Moreover, Figure \ref{MUSCL2_TVD_RK44_s} shows that not only the RK44 solution, but also the $\underline{q}^n + \Delta t \underline{R}^{j=s}$ term maintained the TVD property up to machine precision. Furthermore, Table \ref{table_MUSCL2_TVD} provides $c^s$ and $c^p$ values for tested RK methods, demonstrating again that for all of considered RK methods, whether SSP or non-SSP, both $c^s$ and $c^p$ parameters are higher than zero, as expected from the theory.

\begin{table}[htbp]
\centering
\begin{tabular}{m{3cm} m{1cm} m{1cm} m{1cm}}
\hline
RK scheme & $c^{ssp}$ & $c^s$ & $c^p$ \\ [0.5ex] 
\hline
forward Euler & $1$ &$1.3$ & $1.3$ \\
midpoint rule & $0$ &$1.3$ & $1.7$ \\
SSPRK33 & $1$ &$0.2$ & $1.3$ \\
RK31 & $0$ & $1.3$ & $2.0$ \\
RK44 & $0$ &$1.3$ & $1.7$ \\
\hline
\end{tabular}
\caption{$c^s$ and $c^p$ values for the the Burgers' problem discretized with the second-order MUSCL scheme .}
\label{table_MUSCL2_TVD}
\end{table}

\subsection{Leblanc shocktube}\label{example_Leblanc}

For the compressible Euler equations, the Leblanc shocktube problem is particularly challenging as discretizations without limiting often fail due to the emergence of negative density and pressure. We follow the Leblanc shocktube problem from Lin et al. \cite{linPositivityPreservingStrategy2023a}, where the domain is $[0,1]$, with the initial condition 
\begin{equation}
\underline{q}_0(x)= \left\{ \begin{array}{lcr}
\underline{q}_L, \; & \mbox{for} \; x<x_0, \\
\underline{q}_R, \; & \mbox{for} \; x\geq x_0, \end{array} \right.
\end{equation}

\begin{equation}
\underline{q}_L= \begin{bmatrix}
\rho_L \\
u_L\\
p_L
\end{bmatrix}= \begin{bmatrix}
1.0\\
0.0\\
(\gamma -1)\;0.1
\end{bmatrix}, \quad \underline{q}_R= \begin{bmatrix}
\rho_R \\
u_R\\
p_R
\end{bmatrix}= \begin{bmatrix}
10^{-3}\\
0.0\\
(\gamma -1)\;10^{-10}
\end{bmatrix},
\end{equation}
with $x_0= 0.33$, $\gamma= 5/3$. This test case has an exact solution as provided in \cite{linPositivityPreservingStrategy2023a}. The domain, according to \cite{linPositivityPreservingStrategy2023a}, is uniformly discretized into $NE$ elements with $(N+1)$-point Gauss-Lobatto nodes, where a Discontinuous Galerkin (DG) scheme limited by a subcell local Lax-Friedrichs method becomes positivity preserving with the forward Euler method. 

\begin{figure}[htbp]

\begin{subfigure}[b]{0.5\textwidth}
\centering
\includegraphics[width=\textwidth]{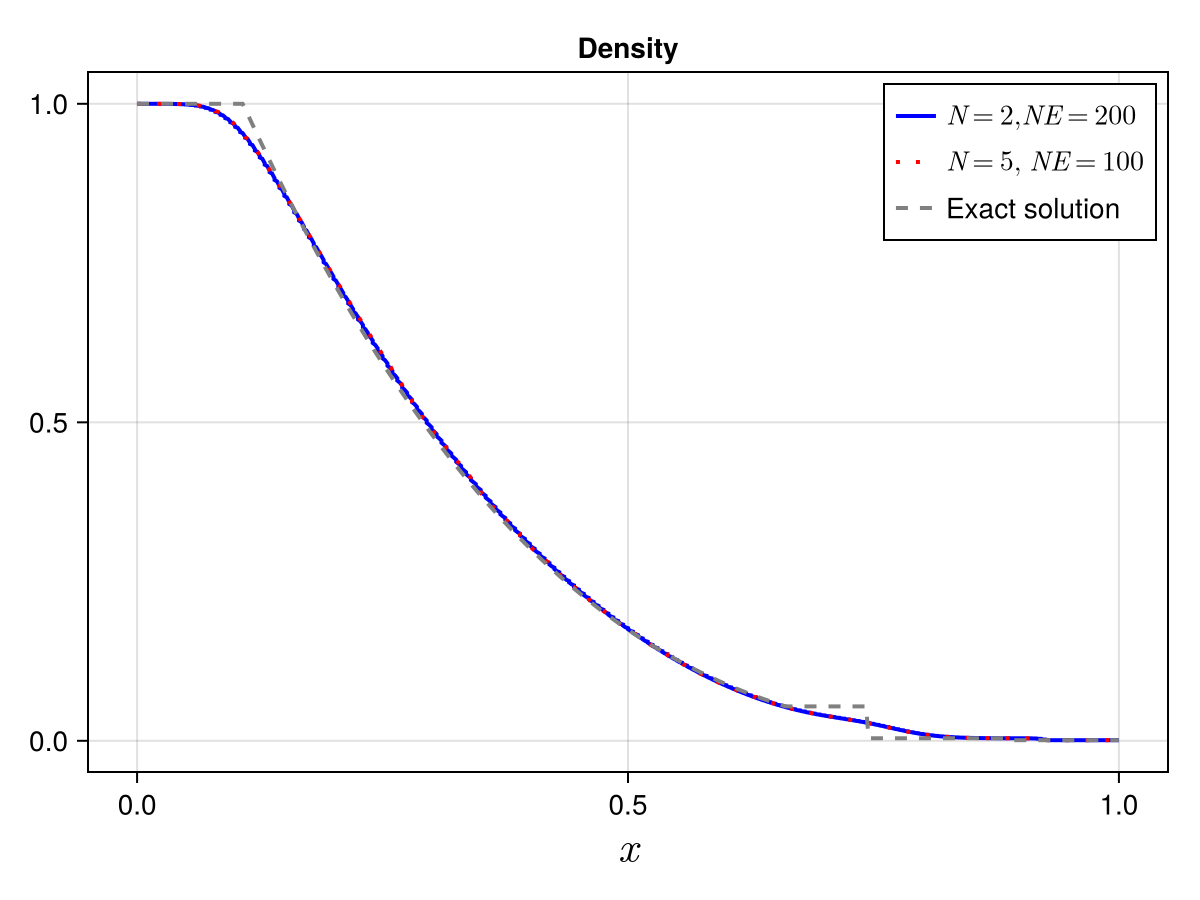}
\caption{Solution}
\end{subfigure}
\begin{subfigure}[b]{0.5\textwidth}
\centering
\includegraphics[width=\textwidth]{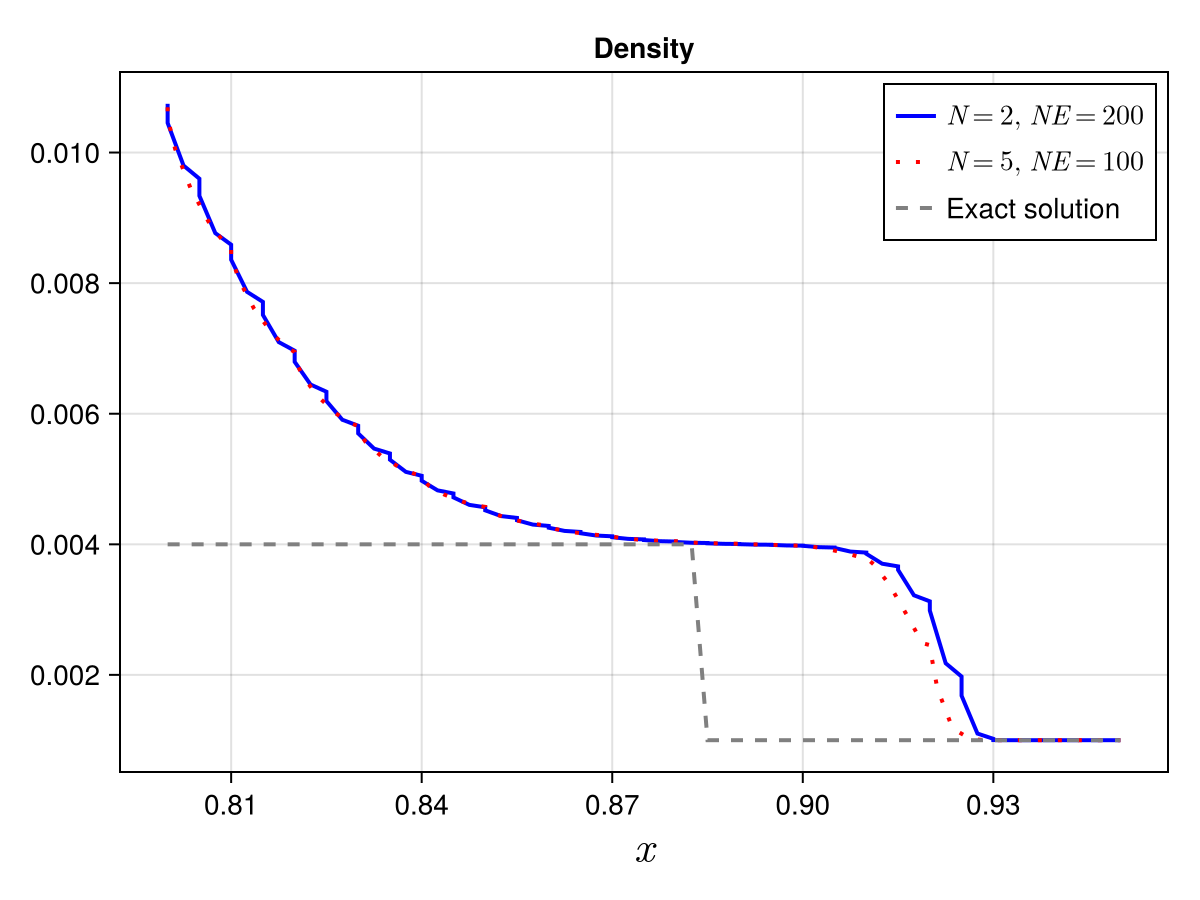}
\caption{Zoomed}
\end{subfigure}

\caption{Density distribution for the Leblanc shocktube problem with the local Lax-Friedrichs discretization, integrated using RK44 with $\Delta t= \Delta t_{FE}$ up to $T= 2/3$ .}
\label{Fig_Leblanc_rho_low}
\end{figure}

\begin{figure}[htbp]

\begin{subfigure}[b]{0.5\textwidth}
\centering
\includegraphics[width=\textwidth]{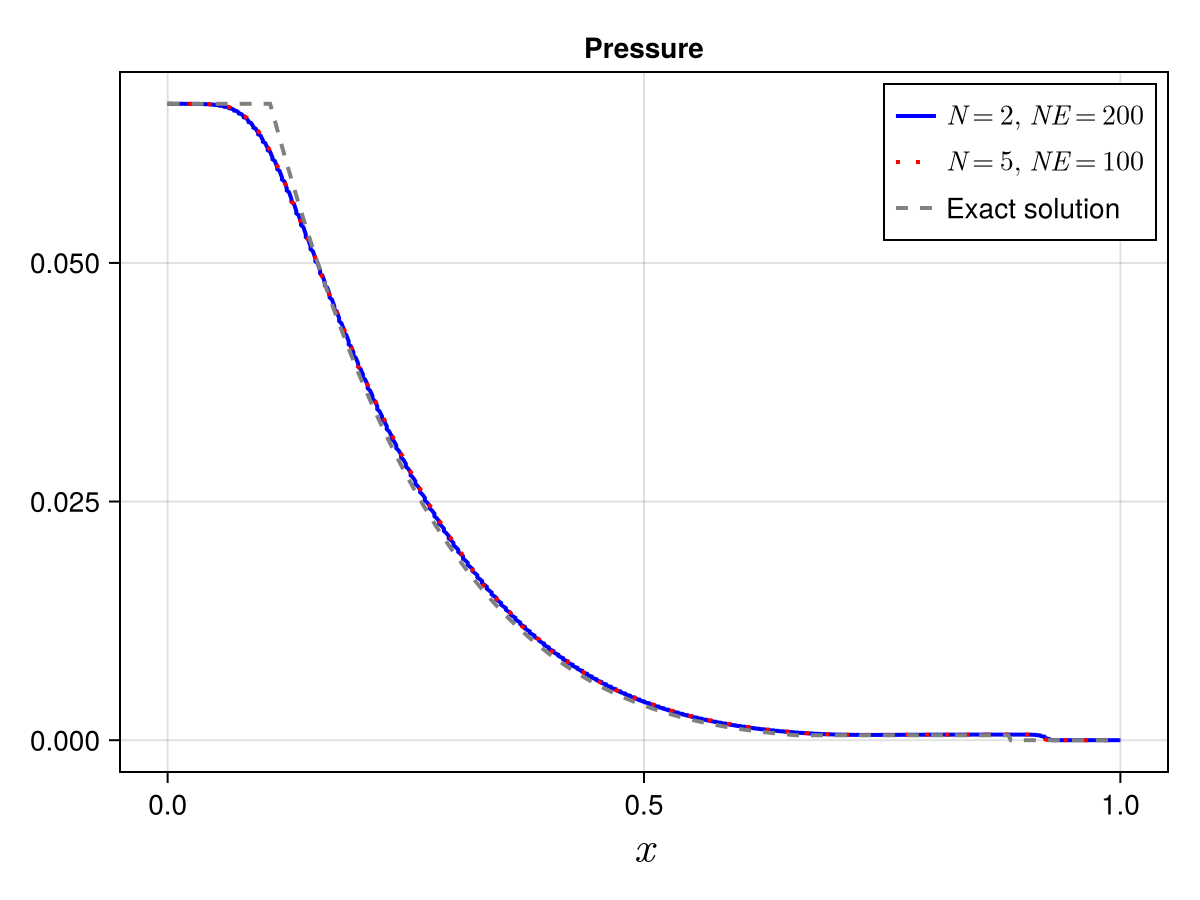}
\caption{Solution}
\end{subfigure}
\begin{subfigure}[b]{0.5\textwidth}
\centering
\includegraphics[width=\textwidth]{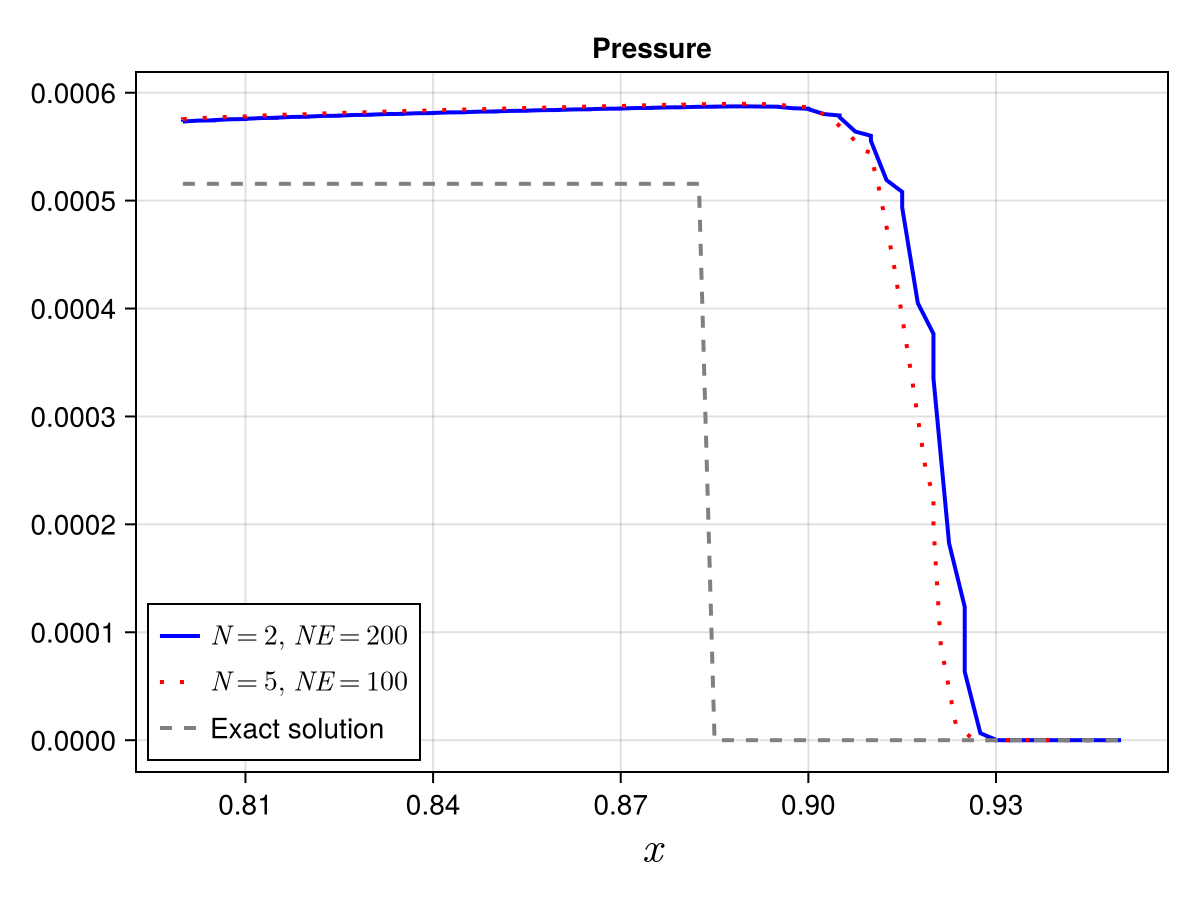}
\caption{Zoomed}
\end{subfigure}

\caption{Pressure distribution for the Leblanc shocktube problem with the local Lax-Friedrichs discretization, integrated using RK44 with $\Delta t= \Delta t_{FE}$ up to $T= 2/3$ .}
\label{Fig_Leblanc_pressure_low}
\end{figure}

We first consider the low-order local Lax-Friedrichs spatial discretization alone. In this case, density and pressure remain positive with the forward Euler method when $\Delta t\leq \Delta t_{FE}$, where $\Delta t_{FE}>0$. For the number of elements $NE=200$ and $NE=100$, corresponding to $N=2$ and $N=5$ respectively, we integrate the problem using RK44 with $\Delta t= \Delta t_{FE}$ up to a final time of $T= 2/3$. As demonstrated in Figures \ref{Fig_Leblanc_rho_low} and \ref{Fig_Leblanc_pressure_low}, both density and pressure remain positive, even though the SSP time-step limit for RK44 is zero. Furthermore, Tables \ref{table_leblanc_200el_low} and \ref{table_leblanc_100el_low} provide $c^s$ and $c^p$ values for all RK methods tested. These results indicate that $\underline{q}^n + \Delta t \underline{R}^j$ terms and non-SSP RK schemes maintained positivity with positive step sizes comparable to the forward Euler limit.

 \begin{table}[htbp]
\centering
\begin{tabular}{m{3cm} m{1cm} m{1cm} m{1cm}}
\hline
RK scheme & $c^{ssp}$ & $c^s$ & $c^p$ \\ [0.5ex] 
\hline
forward Euler & $1$ & $2.7$  &$2.7$ \\
midpoint rule & $0$ & $2.7$ & $2.7$ \\
SSPRK33 & $1$ & $2.0$ & $3.0$ \\
RK31 & $0$ & $2.9$ & $3.0$ \\
RK44 & $0$ & $2.4$ & $2.4$ \\
\hline
\end{tabular}
\caption{Positivity preserving limits for the Leblanc shocktube problem with $N=2$ and $NE=200$, and the local Lax-Friedrichs discretization.}
\label{table_leblanc_200el_low}
\end{table}

 \begin{table}[htbp]
\centering
\begin{tabular}{m{3cm} m{1cm} m{1cm} m{1cm}}
\hline
RK scheme & $c^{ssp}$ & $c^s$ & $c^p$ \\ [0.5ex] 
\hline
forward Euler & $1$ & $2.8$  &$2.8$ \\
midpoint rule & $0$ & $2.7$ & $2.7$ \\
SSPRK33 & $1$ & $2.9$ & $2.9$ \\
RK31 & $0$ & $3.3$ & $3.4$ \\
RK44 & $0$ & $2.4$ & $3.2$ \\
\hline
\end{tabular}
\caption{Positivity preserving limits for the Leblanc shocktube problem with $N=5$ and $NE=100$, and the local Lax-Friedrichs discretization.}
\label{table_leblanc_100el_low}
\end{table}

\begin{figure}[htbp]

\begin{subfigure}[b]{0.4\textwidth}
\centering
\includegraphics[width=\textwidth]{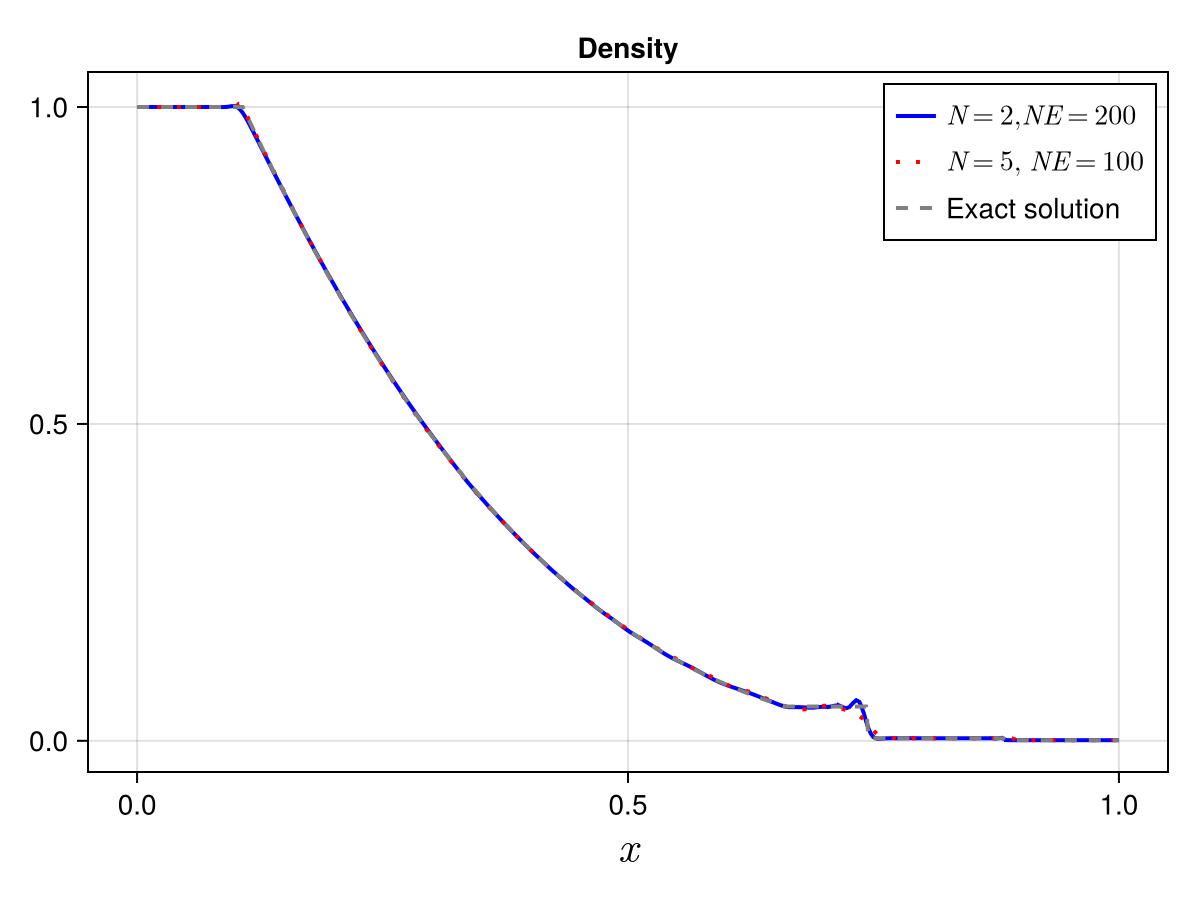}
\caption{Solution}
\end{subfigure}
\begin{subfigure}[b]{0.4\textwidth}
\centering
\includegraphics[width=\textwidth]{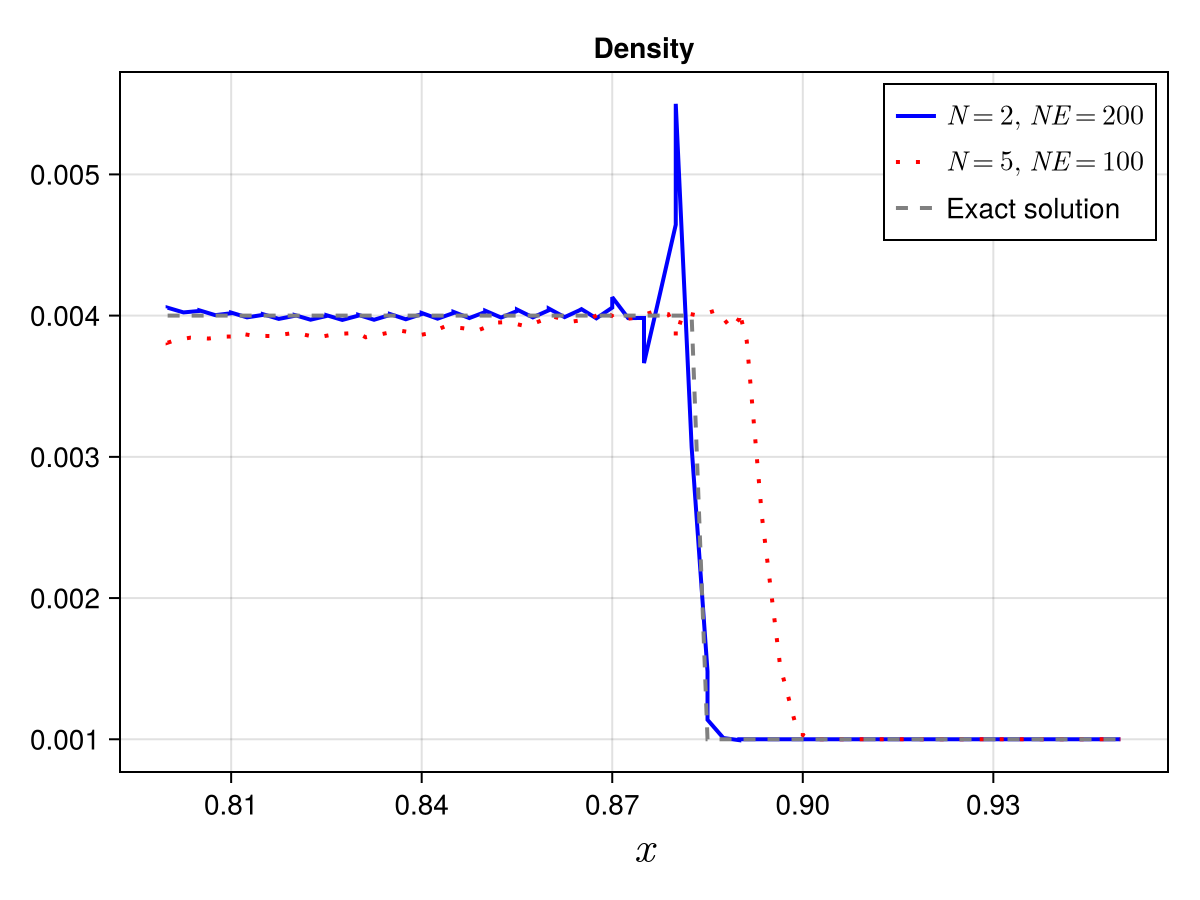}
\caption{Zoomed}
\end{subfigure}

\caption{Density distribution for the Leblanc shocktube problem discretized with a subcell-limited DG scheme and integrated using RK44 with $\Delta t= \Delta t_{FE}$ up to $T= 2/3$. }
\label{Fig_Leblanc_rho_high}
\end{figure}

\begin{figure}[htbp]

\begin{subfigure}[b]{0.4\textwidth}
\centering
\includegraphics[width=\textwidth]{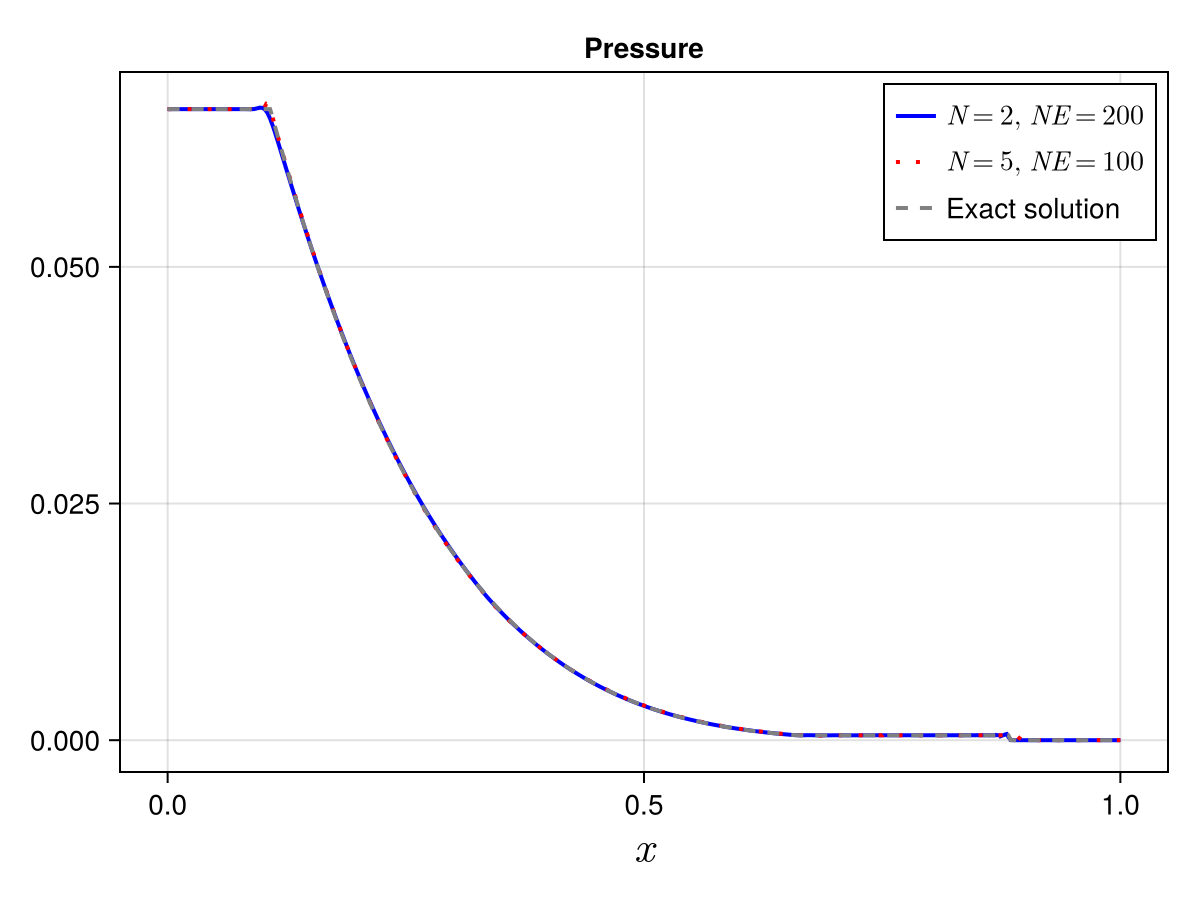}
\caption{Solution}
\end{subfigure}
\begin{subfigure}[b]{0.4\textwidth}
\centering
\includegraphics[width=\textwidth]{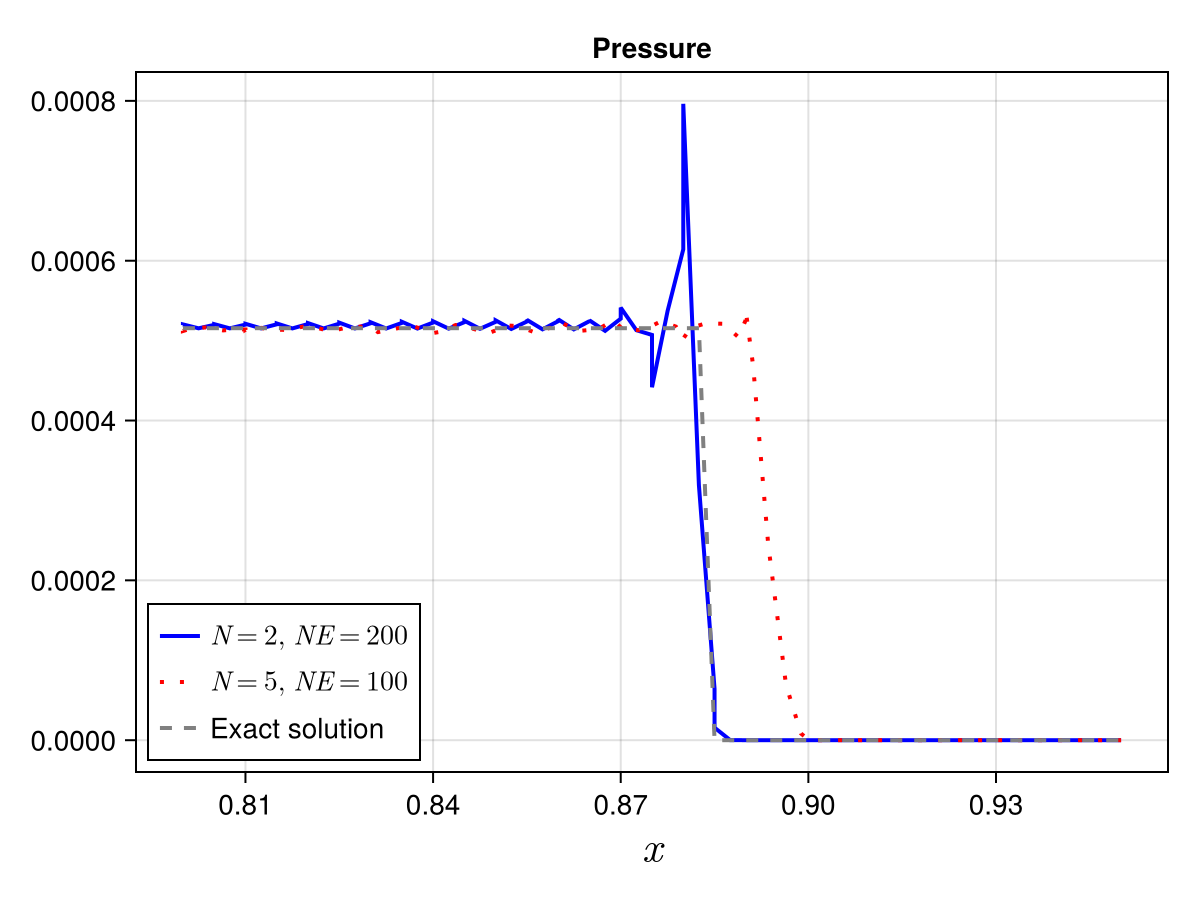}
\caption{Zoomed}
\end{subfigure}

\caption{Pressure distribution for the Leblanc shocktube problem discretized with a subcell-limited DG scheme and integrated using RK44 with $\Delta t= \Delta t_{FE}$ up to $T= 2/3$  .}
\label{Fig_Leblanc_pressure_high}
\end{figure}

Moreover, we test here a Discontinuous Galerkin (DG) scheme limited by a subcell-based Lax-Friedrichs method, according to \cite{linPositivityPreservingStrategy2023a}, to achieve high order spatial accuracy while maintaining positivity of density and pressure with forward Euler time integration. Although we have not provided theoretical proofs for this spatial discretization, we experimentally show that, also with this spatial discretization, positivity is maintained with RK methods satisfying the Assumption \ref{Assumption_RK} when step size is sufficiently small. Figures \ref{Fig_Leblanc_rho_high} and \ref{Fig_Leblanc_pressure_high} present the density and pressure distributions obtained by integrating the problem using RK44 with $\Delta t= \Delta t_{FE}$ up to $T=2/3$. Furthermore, Tables \ref{table_leblanc_200el_high} and \ref{table_leblanc_100el_high} confirm that the terms $\underline{q}^n + \Delta t \underline{R}^j$ and non-SSP schemes again maintained positivity of density and pressure with positive time steps.

\begin{table}[htbp]
\centering
\begin{tabular}{m{3cm} m{1cm} m{1cm} m{1cm}}
\hline
RK scheme & $c^{ssp}$ & $c^s$ & $c^p$ \\ [0.5ex] 
\hline
forward Euler & $1$ & $3.2$  &$3.2$ \\
midpoint rule & $0$ & $2.7$ & $4.1$ \\
SSPRK33 & $1$ & $2.2$ & $3.8$ \\
RK31 & $0$ & $3.2$ & $3.5$ \\
RK44 & $0$ & $1.6$ & $1.6$ \\
\hline
\end{tabular}
\caption{Positivity preserving limits for the Leblanc shocktube problem discretized using a subcell-limited DG with $N=2$ and $NE=200$.}
\label{table_leblanc_200el_high}
\end{table}

\begin{table}[htbp]
\centering
\begin{tabular}{m{3cm} m{1cm} m{1cm} m{1cm}}
\hline
RK scheme & $c^{ssp}$ & $c^s$ & $c^p$ \\ [0.5ex] 
\hline
forward Euler & $1$ & $3.0$ &$3.0$ \\
midpoint rule & $0$ & $3.0$ & $3.0$ \\
SSPRK33 & $1$ & $3.1$ & $3.3$ \\
RK31 & $0$ & $3.4$ & $3.7$ \\
RK44 & $0$ & $2.5$ &$3.3$ \\
\hline
\end{tabular}
\caption{Positivity preserving limits for the Leblanc shocktube problem discretized using a subcell-limited DG with $N=5$ and $NE=100$.}
\label{table_leblanc_100el_high}
\end{table}

\section{Conclusions}

In this work we introduced a novel theoretical framework for analyzing the SSP property of RK schemes. Instead of relying on the traditional Shu-Osher form, we used an alternative RK representation that enables assessing the stability of RK methods that may not be SSP. Using this approach, we theoretically demonstrated that a certain class of RK methods, even if not being SSP, maintains stability in convex functionals with a step size limit greater than zero. Regarding positivity of density and pressure in the Euler equations, we proved that, with the Lax-Friedrichs spatial discretization, these RK schemes preserve positivity for sufficiently small time steps. Moreover, with the first-order upwind and the second-order minmod-based MUSCL schemes, these RK methods maintain the TVD property of Burgers' problem with positive time steps. Numerical experiments confirm the theoretical results and, interestingly, show that the stability time-step limits for the considered non-SSP schemes are comparable, or even larger than, the forward Euler time-step limit.

The significance of this work is that it formally proves the stability time-step limit for a range of non-SSP RK methods, such as RK44, for some well-known nonlinear hyperbolic spatial discretizations is greater than zero. This means that, at least for certain nonlinear hyperbolic problems, it is not necessary to employ SSP schemes and one can employ a larger class of RK methods. While experiments suggest that the stability step-size limits for the tested non-SSP RK methods are not relatively small compared to the forward Euler step size limit, the theoretical framework presented here relies on the use of a sufficiently small time step and does not provide an upper bound for that. A possible direction for the future work is developing a strategy to approximate the stability step-size limits for considered problems. Moreover, although the theoretical results for convex functionals do not rely on specific problems, the analysis of positivity preservation and TVD stability in this work focused on specific spatial discretizations. Another future work is to extend the theoretical findings to other spatial discretizations, such as the subcell-limited DG schemes tested in the Example \ref{example_Leblanc}.

\section*{Acknowledgments}
The authors acknowledge financial support from the Natural Sciences and Engineering Research Council of Canada (NSERC) under grants RGPAS-2017-507988 and RGPIN-2017-06773, the NOVA – FRQNT-NSERC program with the grant number 2023-314433, and Environment and Climate Change Canada. Moreover, Mohammad R. Najafian acknowledges the support from Concordia University through a Graduate Fellowship.

We thank Dr. Sigal Gottlieb for the fruitful discussions and valuable comments. Furthermore, we acknowledge that to reproduce the experiments presented in the paper \cite{linPositivityPreservingStrategy2023a}, the open-access codes provided by this paper are used in this work.
%\input{Data Statement}
%\section*{Data Statement}

%\bibliographystyle{unsrt}
\bibliographystyle{elsarticle-num}

\bibliography{biblio_ssp}

\end{document}